\documentclass{elsart}
\usepackage{amsfonts}
\begin{document}
%
%----------------------START OF LOCAL DEFINITIONS----------------------------
%
\def\myp#1#2{#1}  %% #1=preprint version, #2=submitted version
\def\mypr#1#2{#2} %% #1=equation nr=nr,   #2=eqno=sec.nr
\mypr{\renewcommand{\theequation}{\arabic{equation}}}%
{\renewcommand{\theequation}{\arabic{section}.\arabic{equation}}}
%------------------------------------------------------------------
\myp{\def\s{\ }}{\def\s{}}
\myp{\def\Journal#1#2#3#4#5#6#7%
{{\if#1-{}\else{\textit{#1}\ }\fi}{#2}{\if#3-{}\else{ #3}\fi}%
\ {\bf #4},{ pp.\ #5--#6} (#7)}}
{\def\Journal#1#2#3#4#5#6#7{{\if#3-{#2\ {\bf#4}\ (#7)\ #5}%
\else{#2\ #3\ {\bf#4} (#7) #5}\fi}}}
%% #1=title, #2=journal, #3=series, #4=vol, #5,6=pp, #7=year
%------------------------------------------------------------------
\myp{\def\InBook#1#2#3#4#5#6#7#8%
{\textit{#1}, in \textit{#2},{\if#3-{}\else{ #3,}\fi}%
{\if#4-{}\else{ Vol.\ #4,}\fi} #5, eds.\ (#6), pp.\ #7--#8}}
{\def\InBook#1#2#3#4#5#6#7#8%
{in \textit{#2},{\if#3-{}\else{ #3,}\fi}%
{\if#4-{}\else{ Vol.\ #4,}\fi} #5, eds.\ (#6), p.\ #7}}
%% #1=title, #2=book, #3=series, #4=vol, #5=eds, #6=publ&year, #7,8=pp
%------------------------------------------------------------------
\catcode`\@=11 %% Based on \underset in amstex.sty.
\def\myunderset#1#2{{\begingroup
  \setbox\z@\hbox{\thinmuskip0mu
    \medmuskip\m@ne mu\thickmuskip\@ne mu
    \setbox\tw@\hbox{$#2\m@th$}\kern-\wd\tw@
    ${}#1{}\m@th$}%
  \edef\my@tempa{\endgroup\let\noexpand\relax
    \ifdim\wd\z@<\z@ \mathbin
    \else\ifdim\wd\z@>\z@ \mathrel
    \else \relax\fi\fi}%
  \my@tempa%
  \relax{\mathop{\kern\z@#2}\limits_{#1}}}}
\catcode`\@=12
%------------------------------------------------------------------
\def\halfs{{\scriptstyle\frac{1}{2}}}
\def\half{{\textstyle\frac{1}{2}}}
\def\f{{\rm f}}
\def\wb{{\overline W}}
\def\wp{W\vphantom{\overline W}}
\def\wf{W^{(\f)}}
\def\wbf{{\overline W}\vphantom{W}^{(\f)}}
\def\wu{W\vphantom{\wbf}}
\def\ualpha{{\alpha}^{\vphantom{(\f)}}}
\def\balpha{{\overline\alpha}\vphantom{\alpha}^{\vphantom{(\f)}}}
\def\alphaf{{\alpha}^{(\f)}}
\def\balphaf{{\overline\alpha}\vphantom{\alpha}^{(\f)}}
\def\ubeta{{\beta}\vphantom{\alpha}^{\vphantom{(\f)}}}
\def\bbeta{{\overline\beta}\vphantom{\alpha}^{\vphantom{(\f)}}}
\def\betaf{{\beta}\vphantom{\alpha}^{(\f)}}
\def\bbetaf{{\overline\beta}\vphantom{\alpha}^{(\f)}}
\def\bA{{\bar A}}
\def\Af{A\vphantom{\alpha}^{(\f)}}
\def\bAf{{\bar A}\vphantom{\balpha}^{(\f)}}
\def\bLf{{\bar L}^{(\f)}}
\def\bC{{\overline C}}
\def\bD{{\overline D}}
\def\tb{{\tilde\theta}}
\def\pb{{\tilde\phi}}
\def\bkappa{{\overline\kappa}}
\def\ii{{\rm i}}
\def\B{{\rm B}}
\def\F{{\rm F}}
\def\O{{\rm O}}
\font\sasi=cmssi12
\def\R{\mbox{\sasi R}}
\def\sign{{\rm sign}}
\def\vpWW{\vphantom{$W_{W_{W_W}}$}}
\def\vpOOO{\vphantom{1^{1^1}}}
\def\vpOO{\vphantom{1^1}}
%------------------------------------------------------------------
\def\p{{\rm p}}
\def\q{{\rm q}}
\def\r{{\rm r}}
\def\pq{{\rm pq}}
\def\pr{{\rm pr}}
\def\qp{{\rm qp}}
\def\qr{{\rm qr}}
\def\rp{{\rm rp}}
\def\rq{{\rm rq}}
\def\pqr{{\rm pqr}}
\def\qq{{\rm qq}}
\def\bfp{\hbox{\bf p}}
\def\bfq{\hbox{\bf q}}
%------------------------------------------------------------------
\def\H#1{{}_{#1}{\rm H}_{#1}}
\def\hypd#1#2#3#4{{}_{#1}{\rm H}_{#1}%
\bigg[{#2\atop#3}\bigg|\,#4\,\bigg]}
%------------------------------------------------------------------
%
%              %%%%%%%    %%%%%        %%%%%%    %%%%%   %     %
%              %      %  %             %     %  %     %   %   %
%              %      %  %             %     %  %     %    % %
%              %%%%%%%    %%%%%        %%%%%%   %     %     %
%              %               %       %     %  %     %    % %
%              %               %       %     %  %     %   %   %
%              %         %%%%%%        %%%%%%    %%%%%   %     %
%
%              By Jean Orloff
%              Comments & suggestions by e-mail: ORLOFF@surya11.cern.ch
%              No modification of this file allowed if not e-sent to me.
%
% A simple way to measure the size of encapsulated postscript figures
%   from inside TeX, and to use it for automatically formatting texts
%   with inserted figures. Works both under Plain TeX-based macros
%   (Phyzzx, Harvmac, Psizzl, ...) and LaTeX environment.
% Provides exactly the same result on any PostScript printer provided
%   the single instruction \psfor... is changed to fit the needs of the
%   particular dvi->ps translator used.
% History:
%   1.31: adds \psforDVIALW(?)
%   1.30: adds \splitfile & \joinfiles for multi-file management
%   1.24: fix error handling & add \psonlyboxes
%   1.23: adds \putsp@ce for OzTeX fix
%   1.22: makes \drawingBox \global for use in Phyzzx
%   1.21: accepts %%BoundingBox: (atend)
%   1.20: tries to add \psfordvitps for the TeXPS package.
%   1.10: adds \psforoztex, error handling...
%2345678 1 2345678 2 2345678 3 2345678 4 2345678 5 2345678 6 2345678 7 23456789
%
\def\temp{1.31}
\let\tempp=\relax
\expandafter\ifx\csname psboxversion\endcsname\relax
  \message{version: \temp}
\else
    \ifdim\temp cm>\psboxversion cm
      \message{version: \temp}
    \else
      \message{psbox(\psboxversion) is already loaded: I won't load
        psbox(\temp)!}
      \let\temp=\psboxversion
      \let\tempp= 
    \fi
\fi
\tempp
\let\psboxversion=\temp
\catcode`\@=11
% Every macro likes a little privacy...
%
% Some common defs
%
\def\execute#1{#1}% NOT stupid: cs in #1 are then identified BEFORE execution
\def\psm@keother#1{\catcode`#112\relax}% borrowed from latex
\def\executeinspecs#1{%
\execute{\begingroup\let\do\psm@keother\dospecials\catcode`\^^M=9#1\endgroup}}
%
%Trying to tame the variety of \special commands for Postscript: the
%  universal internal command \PSspeci@l##1##2 takes ##1 to be the
%  filename and ##2 to be the integer scale factor*1000 (as for usual
%   TeX \scale commands)
%
\def\psfortextures{%     For TeXtures on the Macintosh
%-----------------
\def\PSspeci@l##1##2{%
\special{illustration ##1\space scaled ##2}%
}}
\def\psfordvitops{%      For the DVItoPS converter on IBM mainframes
%----------------
\def\PSspeci@l##1##2{%
\special{dvitops: import ##1\space \the\drawingwd \the\drawinght}%
}}
\def\psfordvips{%      For DVIPS converter on VAX, UNIX and PC's
%--------------
\def\PSspeci@l##1##2{%
%    \special{/@scaleunit 1000 def}% never read dox without trying!
\d@my=0.1bp \d@mx=\drawingwd \divide\d@mx by\d@my%
\includegraphics{##1\space}%
}}
\def\psforoztex{%        For the OzTeX shareware on the Macintosh
%--------------
\def\PSspeci@l##1##2{%
\special{##1 \space
      ##2 1000 div dup scale
      \putsp@ce{\number-\psllx} \putsp@ce{\number-\pslly} translate
}%
}}
\def\putsp@ce#1{#1 }
\def\psfordvitps{%       From the UNIX TeXPS package, vers.>3.12
%---------------
% Convert a dimension into the number \psn@sp (in scaled points)
\def\psdimt@n@sp##1{\d@mx=##1\relax\edef\psn@sp{\number\d@mx}}
\def\PSspeci@l##1##2{%
% psfig.psr contains the def of "startTexFig": if you can locate it
% and include the correct pathname, it should work
\special{dvitps: Include0 "psfig.psr"}% contains def of "startTexFig"
\psdimt@n@sp{\drawingwd}
\special{dvitps: Literal "\psn@sp\space"}
\psdimt@n@sp{\drawinght}
\special{dvitps: Literal "\psn@sp\space"}
\psdimt@n@sp{\psllx bp}
\special{dvitps: Literal "\psn@sp\space"}
\psdimt@n@sp{\pslly bp}
\special{dvitps: Literal "\psn@sp\space"}
\psdimt@n@sp{\psurx bp}
\special{dvitps: Literal "\psn@sp\space"}
\psdimt@n@sp{\psury bp}
\special{dvitps: Literal "\psn@sp\space startTexFig\space"}
\special{dvitps: Include1 "##1"}
\special{dvitps: Literal "endTexFig\space"}
}}
\def\psforDVIALW{%   Try for dvialw, a UNIX public domain
%---------------
\def\PSspeci@l##1##2{
\special{language "PS"
literal "##2 1000 div dup scale"
include "##1"}}}
\def\psonlyboxes{%     Draft-like behaviour if none of the others works
%---------------
\def\PSspeci@l##1##2{%
\at(0cm;0cm){\boxit{\vbox to\drawinght
  {\vss
  \hbox to\drawingwd{\at(0cm;0cm){\hbox{(##1)}}\hss}
  }}}
}%
}
\def\psloc@lerr#1{%
\let\savedPSspeci@l=\PSspeci@l%
\def\PSspeci@l##1##2{%
\at(0cm;0cm){\boxit{\vbox to\drawinght
  {\vss
  \hbox to\drawingwd{\at(0cm;0cm){\hbox{(##1) #1}}\hss}
  }}}
\let\PSspeci@l=\savedPSspeci@l% restore normal output for other figs!
}%
}
%
%\def\psfor...  add your own!
%
%  \ReadPSize{PSfilename} reads the dimensions of a PostScript drawing
%      and stores it in \drawinght(wd)
\newread\pst@mpin
\newdimen\drawinght\newdimen\drawingwd
\newdimen\psxoffset\newdimen\psyoffset
\newbox\drawingBox
\newif\ifNotB@undingBox
\newhelp\PShelp{Proceed: you'll have a 5cm square blank box instead of
your graphics (Jean Orloff).}
\def\@mpty{}
\def\s@tsize#1 #2 #3 #4\@ndsize{
  \def\psllx{#1}\def\pslly{#2}%
  \def\psurx{#3}\def\psury{#4}%  needed by a crazyness of dvips!
  \ifx\psurx\@mpty\NotB@undingBoxtrue% this is not a valid one!
  \else
    \drawinght=#4bp\advance\drawinght by-#2bp
    \drawingwd=#3bp\advance\drawingwd by-#1bp
%  !Units related by crazy factors as bp/pt=72.27/72 should be BANNED!
  \fi
  }
\def\sc@nline#1:#2\@ndline{\edef\p@rameter{#1}\edef\v@lue{#2}}
\def\g@bblefirstblank#1#2:{\ifx#1 \else#1\fi#2}
\def\psm@keother#1{\catcode`#112\relax}% borrowed from latex
\def\execute#1{#1}% Seems stupid, but cs are identified BEFORE execution
{\catcode`\%=12
\xdef\B@undingBox{%%BoundingBox}
}   %% is not a true comment in PostScript, even if % is!
\def\ReadPSize#1{
 \edef\PSfilename{#1}
 \openin\pst@mpin=#1\relax
 \ifeof\pst@mpin \errhelp=\PShelp
   \errmessage{I haven't found your postscript file (\PSfilename)}
   \psloc@lerr{was not found}
   \s@tsize 0 0 142 142\@ndsize
   \closein\pst@mpin
 \else
   \immediate\write\psbj@inaux{#1,}
   \loop
     \executeinspecs{\catcode`\ =10\global\read\pst@mpin to\n@xtline}
     \ifeof\pst@mpin
       \errhelp=\PShelp
       \errmessage{(\PSfilename) is not an Encapsulated PostScript File:
           I could not find any \B@undingBox: line.}
       \edef\v@lue{0 0 142 142:}
       \psloc@lerr{is not an EPSFile}
       \NotB@undingBoxfalse
     \else
       \expandafter\sc@nline\n@xtline:\@ndline
       \ifx\p@rameter\B@undingBox\NotB@undingBoxfalse
         \edef\t@mp{%
           \expandafter\g@bblefirstblank\v@lue\space\space\space}
         \expandafter\s@tsize\t@mp\@ndsize
       \else\NotB@undingBoxtrue
       \fi
     \fi
   \ifNotB@undingBox\repeat
   \closein\pst@mpin
 \fi
\message{#1}
}
%
% \psboxto(xdim;ydim){psfilename}: you specify the dimensions and
%    TeX uniformly scales to fit the largest one. If xdim=0pt, the
%    scale is fully determined by ydim and vice versa.
%    Notice: psboxes are a real vboxes; couldn't take hbox otherwise all
%    indentation and all cr's would be interpreted as spaces (hugh!).
%
\newcount\xscale \newcount\yscale \newdimen\pscm\pscm=1cm
\newdimen\d@mx \newdimen\d@my
\let\ps@nnotation=\relax
\def\psboxto(#1;#2)#3{\vbox{
   \ReadPSize{#3}
   \divide\drawingwd by 1000
   \divide\drawinght by 1000
   \d@mx=#1
   \ifdim\d@mx=0pt\xscale=1000
         \else \xscale=\d@mx \divide \xscale by \drawingwd\fi
   \d@my=#2
   \ifdim\d@my=0pt\yscale=1000
         \else \yscale=\d@my \divide \yscale by \drawinght\fi
   \ifnum\yscale=1000
         \else\ifnum\xscale=1000\xscale=\yscale
                    \else\ifnum\yscale<\xscale\xscale=\yscale\fi
              \fi
   \fi
   \divide \psxoffset by 1000\multiply\psxoffset by \xscale
   \divide \psyoffset by 1000\multiply\psyoffset by \xscale
   \global\divide\pscm by 1000
   \global\multiply\pscm by\xscale
   \multiply\drawingwd by\xscale \multiply\drawinght by\xscale
   \ifdim\d@mx=0pt\d@mx=\drawingwd\fi
   \ifdim\d@my=0pt\d@my=\drawinght\fi
   \message{scaled \the\xscale}
 \hbox to\d@mx{\hss\vbox to\d@my{\vss
   \global\setbox\drawingBox=\hbox to 0pt{\kern\psxoffset\vbox to 0pt{
      \kern-\psyoffset
      \PSspeci@l{\PSfilename}{\the\xscale}
      \vss}\hss\ps@nnotation}
   \global\ht\drawingBox=\the\drawinght
   \global\wd\drawingBox=\the\drawingwd
   \baselineskip=0pt
   \copy\drawingBox
 \vss}\hss}
  \global\psxoffset=0pt
  \global\psyoffset=0pt% These are local to one figure
  \global\pscm=1cm
  \global\drawingwd=\drawingwd
  \global\drawinght=\drawinght
}}
%
% \psboxscaled{scalefactor*1000}{PSfilename} allows to bypass the
%   rounding errors of TeX integer divisions for situations where the
%   TeX box should fit the original BoundingBox with a precision better
%   than 1/1000.
%
\def\psboxscaled#1#2{\vbox{
  \ReadPSize{#2}
  \xscale=#1
  \message{scaled \the\xscale}
  \divide\drawingwd by 1000\multiply\drawingwd by\xscale
  \divide\drawinght by 1000\multiply\drawinght by\xscale
  \divide \psxoffset by 1000\multiply\psxoffset by \xscale
  \divide \psyoffset by 1000\multiply\psyoffset by \xscale
  \global\divide\pscm by 1000
  \global\multiply\pscm by\xscale
  \global\setbox\drawingBox=\hbox to 0pt{\kern\psxoffset\vbox to 0pt{
     \kern-\psyoffset
     \PSspeci@l{\PSfilename}{\the\xscale}
     \vss}\hss\ps@nnotation}
  \global\ht\drawingBox=\the\drawinght
  \global\wd\drawingBox=\the\drawingwd
  \baselineskip=0pt
  \copy\drawingBox
  \global\psxoffset=0pt
  \global\psyoffset=0pt% These are local to one figure
  \global\pscm=1cm
  \global\drawingwd=\drawingwd
  \global\drawinght=\drawinght
}}
%
%  \psbox{PSfilename} makes a TeX box having the minimal size to
%      enclose the picture
\def\psbox#1{\psboxscaled{1000}{#1}}
%
%
%  \joinfiles file1, file2, ...n \into joinedfilename .
%     makes one file out of many
%  \splitfile joinedfilename
%     the opposite
%
%\def\execute#1{#1}% NOT stupid: cs in #1 are then identified BEFORE execution
%\def\psm@keother#1{\catcode`#112\relax}% borrowed from latex
%\def\executeinspecs#1{%
%\execute{\begingroup\let\do\psm@keother\dospecials\catcode`\^^M=9#1\endgroup}}
%\newread\pst@mpin
\newif\ifn@teof\n@teoftrue
\newif\ifc@ntrolline
\newif\ifmatch
\newread\j@insplitin
\newwrite\j@insplitout
\newwrite\psbj@inaux
\immediate\openout\psbj@inaux=psbjoin.aux
\immediate\write\psbj@inaux{\string\joinfiles}
\immediate\write\psbj@inaux{\jobname,}
%
% We redefine input to keep track of the various files inputted
%
\immediate\let\oldinput=\input
\def\input#1 {
 \immediate\write\psbj@inaux{#1,}
 \oldinput #1 }
\def\empty{}
\def\setmatchif#1\contains#2{
  \def\match##1#2##2\endmatch{
    \def\tmp{##2}
    \ifx\empty\tmp
      \matchfalse
    \else
      \matchtrue
    \fi}
  \match#1#2\endmatch}
\def\warnopenout#1#2{
 \setmatchif{TrashMe,psbjoin.aux,psbjoin.all}\contains{#2}
 \ifmatch
 \else
   \immediate\openin\pst@mpin=#2
   \ifeof\pst@mpin
     \else
     \errhelp{If the content of this file is so precious to you, abort (ie
press x or e) and rename it before retrying.}
     \errmessage{I'm just about to replace your file named #2}
   \fi
   \immediate\closein\pst@mpin
 \fi
 \message{#2}
 \immediate\openout#1=#2}
%  No comments allowed below: % will have an unusual catcode
{
\catcode`\%=12
\gdef\splitfile#1 {
 \immediate\openin\j@insplitin=#1
 \message{Splitting file #1 into:}
 \warnopenout\j@insplitout{TrashMe}
 \loop
   \ifeof
     \j@insplitin\immediate\closein\j@insplitin\n@teoffalse
   \else
     \n@teoftrue
     \executeinspecs{\global\read\j@insplitin to\spl@tinline\expandafter
       \ch@ckbeginnewfile\spl@tinline%Beginning-Of-File-Named:%\endcheck}
     \ifc@ntrolline
     \else
       \toks0=\expandafter{\spl@tinline}
       \immediate\write\j@insplitout{\the\toks0}
     \fi
   \fi
 \ifn@teof\repeat
 \immediate\closeout\j@insplitout}
\gdef\ch@ckbeginnewfile#1%Beginning-Of-File-Named:#2%#3\endcheck{
 \def\t@mp{#1}
 \ifx\empty\t@mp
   \def\t@mp{#3}
   \ifx\empty\t@mp
     \global\c@ntrollinefalse
   \else
     \immediate\closeout\j@insplitout
     \warnopenout\j@insplitout{#2}
     \global\c@ntrollinetrue
   \fi
 \else
   \global\c@ntrollinefalse
 \fi}
\gdef\joinfiles#1\into#2 {
 \message{Joining following files into}
 \warnopenout\j@insplitout{#2}
 \message{:}
 {
 \edef\w@##1{\immediate\write\j@insplitout{##1}}
 \w@{% This text was produced with psbox's \string\joinfiles.}
 \w@{% To decompose and tex it:}
 \w@{%-save this with a filename CONTAINING ONLY LETTERS, and no extensions}
 \w@{% (say, JOINTFIL), in some uncrowded directory;}
 \w@{%-make sure you can \string\input\space psbox.tex (version>=1.3);}
 \w@{%-tex JOINTFIL using Plain, or LaTeX, or whatever is needed by}
 \w@{% the first part in the joining (after splitting JOINTFIL into}
 \w@{% it's constituents, TeX will try to process it as it stands).}
 \w@{\string\input\space psbox.tex}
 \w@{\string\splitfile{\string\jobname}}
 }
 \tre@tfilelist#1, \endtre@t
 \immediate\closeout\j@insplitout}
\gdef\tre@tfilelist#1, #2\endtre@t{
 \def\t@mp{#1}
 \ifx\empty\t@mp
   \else
   \llj@in{#1}
   \tre@tfilelist#2, \endtre@t
 \fi}
\gdef\llj@in#1{
 \immediate\openin\j@insplitin=#1
 \ifeof\j@insplitin
   \errmessage{I couldn't find file #1.}
   \else
   \message{#1}
   \toks0={%Beginning-Of-File-Named:#1}
   \immediate\write\j@insplitout{\the\toks0}
   \executeinspecs{\global\read\j@insplitin to\oldj@ininline}
   \loop
     \ifeof\j@insplitin\immediate\closein\j@insplitin\n@teoffalse
       \else\n@teoftrue
       \executeinspecs{\global\read\j@insplitin to\j@ininline}
       \toks0=\expandafter{\oldj@ininline}
       \let\oldj@ininline=\j@ininline
       \immediate\write\j@insplitout{\the\toks0}
     \fi
   \ifn@teof
   \repeat
   \immediate\closein\j@insplitin
 \fi}
}
% To be put at the end of a file, for making an tar-like file containing
%   everything it used.
\def\autojoin{
 \immediate\write\psbj@inaux{\string\into\space psbjoin.all}
 \immediate\closeout\psbj@inaux
 \input psbjoin.aux
}
%
%  Annotations & Captions etc...
%
%
% \centinsert{anybox} is just a centered \midinsert, but is included as
%    people barely use the original inserts from TeX.
%
\def\centinsert#1{\midinsert\line{\hss#1\hss}\endinsert}
\def\psannotate#1#2{\def\ps@nnotation{#2\global\let\ps@nnotation=\relax}#1}
\def\pscaption#1#2{\vbox{
   \setbox\drawingBox=#1
   \copy\drawingBox
   \vskip\baselineskip
   \vbox{\hsize=\wd\drawingBox\setbox0=\hbox{#2}
     \ifdim\wd0>\hsize
       \noindent\unhbox0\tolerance=5000
    \else\centerline{\box0}
    \fi
}}}
% for compatibility with older versions
\def\psfig#1#2#3{\pscaption{\psannotate{#1}{#2}}{#3}}
\def\psfigurebox#1#2#3{\pscaption{\psannotate{\psbox{#1}}{#2}}{#3}}
%
% \at(#1;#2)#3 puts #3 at #1-higher and #2-right of the current
%    position without moving it (to be used in annotations).
\def\at(#1;#2)#3{\setbox0=\hbox{#3}\ht0=0pt\dp0=0pt
  \rlap{\kern#1\vbox to0pt{\kern-#2\box0\vss}}}
%
% \gridfill(ht;wd) makes a 1cm*1cm grid of ht by wd whose lower-left
%   corner is the current point
\newdimen\gridht \newdimen\gridwd
\def\gridfill(#1;#2){
  \setbox0=\hbox to 1\pscm
  {\vrule height1\pscm width.4pt\leaders\hrule\hfill}
  \gridht=#1
  \divide\gridht by \ht0
  \multiply\gridht by \ht0
  \gridwd=#2
  \divide\gridwd by \wd0
  \multiply\gridwd by \wd0
  \advance \gridwd by \wd0
  \vbox to \gridht{\leaders\hbox to\gridwd{\leaders\box0\hfill}\vfill}}
%
% Useful to measure where to put annotations
\def\fillinggrid{\at(0cm;0cm){\vbox{
  \gridfill(\drawinght;\drawingwd)}}}
%
% \textleftof\anybox: Sample text\endtext
%   inserts "Sample text" on the left of \anybox ie \vbox, \psbox.
%   \textrightof is the symmetric (not documented, too uggly)
% Welcome any suggestion about clean wraparound macros from
%   TeXhackers reading this
%
\def\textleftof#1:{
  \setbox1=#1
  \setbox0=\vbox\bgroup
    \advance\hsize by -\wd1 \advance\hsize by -2em}
\def\textrightof#1:{
  \setbox0=#1
  \setbox1=\vbox\bgroup
    \advance\hsize by -\wd0 \advance\hsize by -2em}
\def\endtext{
  \egroup
  \hbox to \hsize{\valign{\vfil##\vfil\cr%
\box0\cr%
\noalign{\hss}\box1\cr}}}
%
% \frameit{\thick}{\skip}{\anybox}
%    draws with thickness \thick a box around \anybox, leaving \skip of
%    blank around it. eg \frameit{0.5pt}{1pt}{\hbox{hello}}
% \boxit{\anybox} is a shortcut.
\def\frameit#1#2#3{\hbox{\vrule width#1\vbox{
  \hrule height#1\vskip#2\hbox{\hskip#2\vbox{#3}\hskip#2}%
        \vskip#2\hrule height#1}\vrule width#1}}
\def\boxit#1{\frameit{0.4pt}{0pt}{#1}}
\catcode`\@=12 % cs containing @ are unreachable
%
% CUSTOMIZE YOUR DEFAULT DRIVER:
%    Uncomment the line corresponding to your TeX system:
%\psfortextures%     For TeXtures on the Macintosh
%\psforoztex   %     For OzTeX shareware on the Macintosh
%\psfordvitops %     For the DVItoPS converter for TeX on IBM mainframes
 \psfordvips   %     For DVIPS converter on VAX and UNIX
%\psfordvitps  %     For dvitps from TeXPS package under UNIX
%\psforDVIALW  %     For DVIALW, UNIX public domain
%\psonlyboxes  %     Blank Boxes (when all else fails).

%
%----------------------END OF LOCAL DEFINITIONS------------------------------
%
\myp{math.QA/9906029\hfill ITFA-98-15}{}
\runauthor{Helen Au-Yang and Jacques H.H.\ Perk}
\begin{frontmatter}
\title{The Large {\em N} Limits of the Chiral Potts Model}
\author{Helen Au-Yang} \and
\author{Jacques H.H.\ Perk}
\address{Department of Physics, Oklahoma State University, %\\
145 Physical Sciences, Stillwater, OK 74078-3072, USA\thanksref{perm}}
\thanks[perm]{Permanent address. E-mail address: perk@okstate.edu$\,$.}
\address{Department of Mathematics and Statistics,
The University of Melbourne, %\\
Parkville, Victoria 3052, Australia}
%\address{\and}
\address{Institute for Theoretical Physics, University of Amsterdam, %\\
Valckenierstraat 65, 1018 XE Amsterdam, The Netherlands.}
\begin{abstract}
In this paper we study the large-$N$ limits of the integrable
$N$-state chiral Potts model. Three chiral solutions of
the star-triangle equations are derived, with states taken from all
integers, or from a finite or infinite real interval. These solutions
are expected to be chiral-field lattice deformations of parafermionic
conformal field theories. A new two-sided hypergeometric identity is
derived as a corollary.
\end{abstract}
\end{frontmatter}

%%%%%%%%%%%%%%%%%%%%%%%%%%%%%%%%%%%%%%%%%%%%%%%%%%%%%%%%%%%%%%%%%%%%%%%%%%%
%  Section 1
%%%%%%%%%%%%%%%%%%%%%%%%%%%%%%%%%%%%%%%%%%%%%%%%%%%%%%%%%%%%%%%%%%%%%%%%%%%
\section{Introduction\label{sec-intro}}

When the integrable $N$-state chiral Potts model was introduced,
it was the first example of an exactly solvable lattice
model whose Boltzmann weights both require the use of higher-genus
algebraic functions for their uniformization and do not have ``the
difference property" \cite{AMPTY,BPTS,P-th,AMPT,BPA,AP-Ta}. Since then,
much has been written about many aspects of this model and we refer
the reader to the recent review \cite{AP-mf} for more information.
In this paper we shall concentrate our attention on just one aspect,
namely the large $N$ limit. We have written about this once before
\cite{AP-inf}, but we can now present a much more complete and
improved version containing several new results in addition.

The chiral Potts model is a spin model on a two-dimensional lattice
(or more generally a planar graph). At each site (or vertex) of the
lattice (or graph), there is a state variable or ``spin" that takes
on $N$ values $a,b,\cdots=1,2,\cdots,N$, (mod $N$). The Boltzmann
weights are associated with pair interactions along edges. We assume
that there are two types of such weights $W$ and $\wb$, which on a
square lattice would correspond to horizontal and vertical
interactions. We assume also that the weights only depend on
the difference modulo $N$ of the two spin states $a$ and $b$ at
the two endpoints of each edge, which is the Potts property. The
chiral character\footnote{This chiral aspect allows us to mimic
the effect of further-neighbor interactions within the context
of a nearest-neighbor interaction model, see \cite{AP-mf} and
references quoted there.} (handedness or breakdown of parity)
is expressed by ${W(a-b)\not\equiv W(b-a)}$ and can only occur
if ${N>2}$.

The integrable chiral Potts model \cite{AMPTY,BPTS,P-th,AMPT,BPA,AP-Ta}
is a nontrivial generalization of the critical Fateev--Zamolodchikov
model \cite{FZ}. The fact that its rapidity variables lie on a
higher-genus curve makes this model special among the many
solvable lattice models. In spite of this several results exist
for it. Therefore, its large-$N$ limits should provide interesting
generalizations of certain nonchiral $\infty$-state models of Fateev
and Zamolodchikov \cite{FZ,Za,Symanzik}, very different from the
SOS-model of Baxter \cite{Ba-SOS} and the few other $\infty$-state
models \cite{FZ-rot,Gaudin,SU,Shibukawa} that have been introduced.
From existing thermodynamic results for the finite-$N$ case, we
can infer corresponding results for the $N=\infty$ cases that
may be of interest in later studies. We expect, for example, a
direct relation with new integrable chiral-field deformations
of parafermionic conformal field theories.

This paper is organized as follows. In section \ref{sec-cpm} we present
the Boltzmann weights of the integrable chiral Potts model and
its dual model, also adding new details not given in \cite{AP-inf}.
In section \ref{sec-inf} we give the three different large-$N$ limits
of the weights, while treating the more technical details in
Appendix \ref{app-N}. The three corresponding large-$N$ limits of
the star-triangle equations are given in detail in section \ref{sec-lim}.
In section \ref{sec-hyp}, it is shown that the results of the previous
section \ref{sec-lim} imply a new two-sided hypergeometric summation
formula. Finally, a short discussion is given in section \ref{sec-dis}.

%%%%%%%%%%%%%%%%%%%%%%%%%%%%%%%%%%%%%%%%%%%%%%%%%%%%%%%%%%%%%%%%%%%%%%%%%%%
%  Section 2
%%%%%%%%%%%%%%%%%%%%%%%%%%%%%%%%%%%%%%%%%%%%%%%%%%%%%%%%%%%%%%%%%%%%%%%%%%%
\mypr{}{\setcounter{equation}{0}}
\section{Integrable $N$-state Chiral Potts Model\label{sec-cpm}}

In this section we shall review earlier results on our
higher-genus solution of the star-triangle equations for the
chiral Potts model \cite{AMPTY,BPTS,P-th,AMPT,BPA,AP-Ta} and
present in more detail a reparametrization \cite{AP-inf,Ba}
that is particularly suitable for the large-$N$ limit.

%%%%%%%%%%%%%%%%%%%%%%%%%%%%%%%%%%%%%%%%%%%%%%%%%%%%%%%%%%%%%%%%%%%%%%%%%%%
\subsection{Star-Triangle Equation for
Chiral Potts Model\label{susec-ste}}

The $N$-state chiral Potts model can be defined on a general
graph with spin states $a,b,\cdots$ taking values $1,\cdots,N$
on the vertices and Boltzmann weights $W(a,b)=W(a-b)$ associated
with edges. $W(n)$ is periodic in $n$ mod $N$.

In the integrable model, one assumes that there are oriented
straight lines (the rapidity lines) on the medial graph, which
are dashed lines shown in Fig.~\ref{fig-latt} for the case of
a square lattice. They are obtained by connecting the middles
of all pairs of edges (solid lines in the figure) that are
incident to a single site and share a common face. No more
than two rapidity lines meet at any given point. These lines
carry variables $\p,\q,\cdots$ and arrows specifying their
orientations. In nearly all solvable models the weights depend
on the differences of these rapidity variables. 
\begin{figure}[tbp]
\vspace{1pc}
%\begin{center}\mbox{\psboxto(\hsize;0cm){fig1.eps}}\end{center}
\begin{center}\mbox{\psboxto(3in;0cm){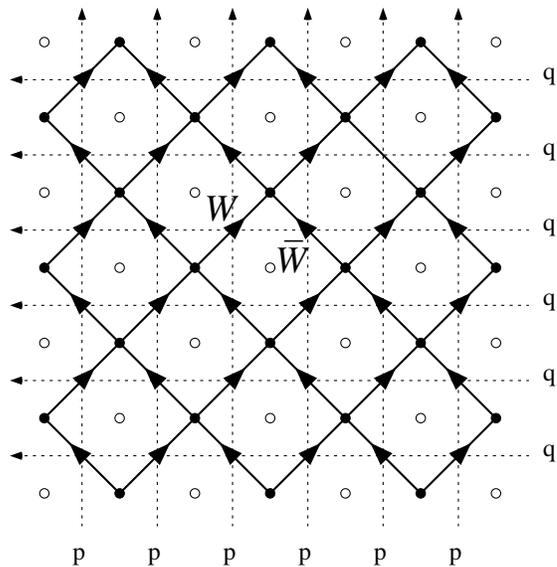}}\end{center}
\vspace{1pc}
\caption{ The square lattice represented by solid lines is diagonally
oriented here, with spins taking $N$ different values at its vertices
and Boltzmann weights $W$ and $\wb$ associated with pair interactions
along the two types of edges. The positions of the spins of the dual
lattice are indicated by open circles. The medial graph is represented
by the dashed oriented horizontal and vertical lines, also called
``rapidity lines," carrying the spectral or rapidity variables
$\p$ and $\q$.}\label{fig-latt}
\end{figure}
For our class of integrable spin-pair interaction models
the weights can be graphically represented as in Fig.~\ref{fig-bw}.
\begin{figure}[tbhp]
\vspace{1pc}
\begin{center}\mbox{\psboxto(3.4in;0cm){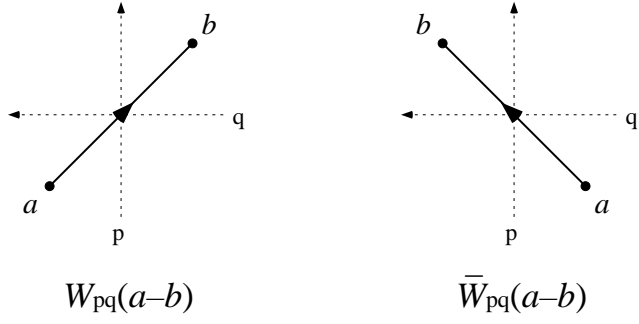}}\end{center}
\vspace{1pc}
\caption{ Boltzmann weights $W_{\pq}(a-b)$ and $\wb_{\pq}(a-b)$.
We need to put an arrow on each edge to distinguish $W_{\pq}(a-b)$
from $W_{\pq}(b-a)$. Note the relative orientation of this arrow
with respect to the orientations of the two rapidity lines in each
case.}\label{fig-bw}
\end{figure}
These weights must satisfy the star-triangle equation
\begin{eqnarray}
&&\sum^{N}_{d=1}\,\wb_{\qr}(b-d)\,W_{\pr}(a-d)\,\wb_{\pq}(d-c)
\nonumber\\
&&\qquad=R_{\pqr}\,W_{\pq}(a-b)\,\wb_{\pr}(b-c)\,W_{\qr}(a-c).
\label{STE}
\end{eqnarray}
Here the factor $R_{\pqr}$ can be determined as \cite{BPA,AP-Ta,MS}
\begin{equation}
R_{\pqr}={F_{\pq}F_{\qr}\over F_{\pr}},\quad
F_{\pq}=\left\{{\prod_{l=1}^{N}\,\sum^{N}_{j=1}\,
\omega^{-jl}\,\wb_{\pq}(j)\over
\prod_{l=1}^{N}\,W_{\pq}(l)}\right\}^{1/N},
\label{STER}
\end{equation}
with
\begin{equation}
\omega\equiv\e^{2\pi{\ii}/N}\equiv\e^{2\pi\sqrt{-1}/N}.
\label{omega}
\end{equation}
The easiest way to derive (\ref{STER}) is to set $a=0$ in (\ref{STE})
and then to take the determinant with respect to the matrix indices $b$
and $c$, leading to determinants of products of diagonal and cyclic
matrices; this argument first appeared in print in \cite{MS}.
The star-triangle equation (\ref{STE}) can be symbolically represented
as in Fig.~\ref{fig-ste}.
\begin{figure}[tbhp]
\vspace{1pc}
%\begin{center}\mbox{\psboxto(\hsize;0cm){fig3.eps}}\end{center}
\begin{center}\mbox{\psboxto(4.5in;0cm){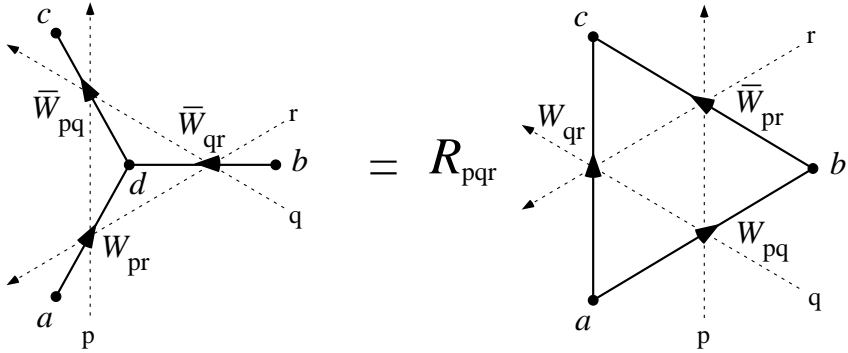}}\end{center}
\vspace{1pc}
\caption{ The star-triangle relations, which allow one to move a
rapidity line $\p$ through a vertex, which is the intersection of
two other rapidity lines $\q$ and $\r$.}\label{fig-ste}
\end{figure}

%%%%%%%%%%%%%%%%%%%%%%%%%%%%%%%%%%%%%%%%%%%%%%%%%%%%%%%%%%%%%%%%%%%%%%%%%%%
\subsection{Weights of Integrable
Chiral Potts Model\label{susec-wcp}}

In \cite{AMPTY,BPTS,P-th,AMPT,BPA,AP-Ta} one family of integrable
chiral Potts models, with an arbitrary number of states per site
$N\ge2$, has been deduced with weights $W_{\pq}(a-b)$ and
$\wb_{\pq}(a-b)$ satisfying the star-triangle equation (\ref{STE})
for all $a,b,c=1,\ldots,N$. These weights are given by
\begin{eqnarray}
&&\frac{W_{\pq}(n)}{ W_{\pq}(0)}=
{\biggl(\frac{\mu_{\p}}{\mu_{\q}}\biggr)}^{\!n}
\prod^{n}_{j=1}\frac{y_{\q}-x_{\p}\omega^j}
{y_{\p}-x_{\q}\omega^j},\nonumber\\
&&\frac{\wb_{\pq}(n)}{\wb_{\pq}(0)}={(\mu_{\p}\mu_{\q})}^n
\prod^{n}_{j=1}\frac{\omega x_{\p}-x_{\q}\omega^j}
{y_{\q}-y_{\p}\omega^j}.
\label{intro1}
\end{eqnarray}
Here,\footnote{To connect with the original homogeneous
notation \cite{BPA,AP-Ta}, we must set
$x_{\p}\equiv a_{\p}/d_{\p}$, $y_{\p}\equiv b_{\p}/c_{\p}$,
$\mu_{\p}\equiv d_{\p}/c_{\p}$, and similarly with ${\p}$ replaced
by ${\q},{\r},\cdots$. A proof that the star-triangle equation
(\ref{STE}) is satisfied is given in the appendix of \cite{AP-Ta}.}
the parameters $\bfp\!\equiv\!(x_{\p},y_{\p},\mu_{\p})$
and\vpWW\
$\bfq\!\equiv\!(x_{\q},y_{\q},\mu_{\q})$ are restricted
by the two periodicity requirements
$W_{\pq}(N+n)=W_{\pq}(n)$, $\wb_{\pq}(N+n)=\wb_{\pq}(n)$, yielding
\begin{equation}
{\biggl(\frac{\mu_{\p}}{\mu_{\q}}\biggr)}^{\!N}=
\,\frac{y_{\p}^N-x_{\q}^N}{y_{\q}^N-x_{\p}^N},
\quad(\mu_{\p}\mu_{\q})^N=
\,\frac{y_{\q}^N-y_{\p}^N}{x_{\p}^N-x_{\q}^N},
\label{intro2}
\end{equation}
which can be recombined as
\begin{equation}
\frac{\mu_{\p}^N x_{\p}^N\pm y_{\p}^N}{1\pm\mu_{\p}^N}=
\frac{\mu_{\q}^N x_{\q}^N\pm y_{\q}^N}{1\pm\mu_{\q}^N}\equiv
\lambda_{\pm},
\label{intro3}
\end{equation}
independent of $\bfp$ and $\bfq$. We write
$\lambda_{\pm}=\pm c(1\mp k')/k$, with $c$ a constant that can be
absorbed by a trivial rescaling of all the
$x_{\p},y_{\p},x_{\q},y_{\q}$ by a common factor that will drop
out of (\ref{intro1}), and with
$k$ and $k'$ numbers related by ${k\vphantom{'}}^2+{k'}^2=1$.
Then the conditions (\ref{intro2}) reduce to
\begin{equation}
\mu_{\p}^N=k'/(1-k\,x_{\p}^N)=(1-k\,y_{\p}^N)/k',\quad
x_{\p}^N+y_{\p}^N=k(1+x_{\p}^N y_{\p}^N).
\label{intro4}
\end{equation}
These equations describe a complex curve, which is the
intersection of two ``Fermat cylinders," and the genus of
this curve is $g=N^2(N-2)+1$. For each Boltzmann weight
the two line (or rapidity) variables $\bfp$ and
$\bfq$ are two points on this higher-genus algebraic
curve, so that the usual difference-variable transformation
cannot be carried out, except for special subcases where
the genus degenerates to $g\le1$. Here the substitutions
$k=0$, $k'=\pm1$ reduce the curve (\ref{intro4}) to a product
of genus-zero curves and the weights degenerate to those of
the self-dual Fateev--Zamolodchikov model \cite{FZ}.

In general the model is not self-dual and the dual weights
are found by Fourier transform. We first note that
the periodic weights in (\ref{intro1}) are of the form
\begin{equation}
{W(n)\over W(0)}=\prod^{n}_{j=1}{x_1-x_2\omega^j 
\over x_4-x_3\omega^j},\quad
{x_1^N-x_2^N\over x_4^N-x_3^N}=1,
\label{intro5}
\end{equation}
so that the linear recursion relation
\begin{equation}
(x_4-x_3\omega^n)W(n)=(x_1-x_2\omega^n)W(n-1)
\label{intro6}
\end{equation}
is satisfied. We can apply the Fourier (or duality) transformation%
\footnote{Here we added a normalization factor $N^{-1}$ that is needed
in the $N\to\infty$ limit.}
\begin{equation}
\wf(m)\equiv N^{-1}\sum_{n=0}^{N-1}\omega^{-mn}W(n),
\label{intro7}
\end{equation}
leading to
\begin{eqnarray}
&&x_4\wf(m)-x_3\wf(m-1)\nonumber\\
&&\qquad=\omega^{-m}x_1\wf(m)-\omega^{1-m}x_2\wf(m-1),
\label{intro8}
\end{eqnarray}
or
\begin{equation}
{\wf(n)\over\wf(0)}=\prod^{n}_{j=1}
{\omega x_2-x_3\omega^j\over x_1-x_4\omega^j}.
\label{intro9}
\end{equation}
Therefore, the weights dual to (\ref{intro1}) are
\begin{eqnarray}
&&\frac{\wf_{\pq}(n)}{\wf_{\pq}(0)}=\prod^{n}_{j=1}
\frac{\omega\mu_{\p}x_{\p}-\mu_{\q}x_{\q}\omega^j}
{\mu_{\p}y_{\q}-\mu_{\q}y_{\p}\omega^j},\nonumber\\
&&\frac{\wbf_{\pq}(n)}{\wbf_{\pq}(0)}=\prod^{n}_{j=1}
\frac{\omega\mu_{\p}\mu_{\q}x_{\q}-y_{\p}\omega^j}
{\omega\mu_{\p}\mu_{\q}x_{\p}-y_{\q}\omega^j},
\label{intro10}
\end{eqnarray}
which are again both of the form (\ref{intro5}).

%%%%%%%%%%%%%%%%%%%%%%%%%%%%%%%%%%%%%%%%%%%%%%%%%%%%%%%%%%%%%%%%%%%%%%%%%%%
\subsection{Reparametrization\label{susec-repa}}

In order to proceed, we introduce new parameters to describe the
higher-genus curve of rapidity variables. These parameters are real
when the Boltzmann weights $W_{\pq}(a-b)$ and $\wb_{\pq}(a-b)$
are real and positive. We begin with the substitutions\footnote{Our
definitions of $\theta_{\p}$ and $\phi_{\p}$ differ by a factor
$N$ from Baxter's \cite{Ba}. This change of normalization will be
necessary in the large $N$ limit.} 
\begin{equation}
x_{\p}=\e^{\ii\phi_{\p}/N},\quad
y_{\p}=\omega^{\halfs}\e^{\ii\theta_{\p}/N},\quad
x_{\q}=\e^{\ii\phi_{\q}/N},\quad
y_{\q}=\omega^{\halfs}\e^{\ii\theta_{\q}/N},
\label{repa0}
\end{equation}
so that from the last identity in (\ref{intro4}) we find
\begin{equation}
k={\sin\half(\theta_{\p}-\phi_{\p})\over
\sin\half(\theta_{\p}+\phi_{\p})}
={\sin\half(\theta_{\q}-\phi_{\q})\over
\sin\half(\theta_{\q}+\phi_{\q})}.
\label{repa1}
\end{equation}
This is equivalent to
\begin{equation}
\e^{\ii\phi_{\p}}={\e^{\ii\theta_{\p}}+k
\over1+k\e^{\ii\theta_{\p}}},
\label{repa2}
\end{equation}
and similarly with ${\p}$ replaced by ${\q}$. From (\ref{repa2})
we have
\begin{equation}
\cos\phi_{\p}={2k+(1+k^2)\cos\theta_{\p}\over
1+k^2+2k\cos\theta_{\p}},
\quad
\sin\phi_{\p}={(1-k^2)\sin\theta_{\p}\over
1+k^2+2k\cos\theta_{\p}}.
\label{repa3}
\end{equation}
We will also need two parameters that will describe the dual model,
see e.g.\ (\ref{repa15}). The first one is given by
\begin{equation}
\lambda_{\p}\equiv{\theta_{\p}+\phi_{\p}\over2\pi}=
{1\over\pi}\arctan{\sin\theta_{\p}\over\cos\theta_{\p}+k},
\label{repa4}
\end{equation}
where the last step follows from (\ref{repa2}). Also using
(\ref{repa1}) we find
\begin{equation}
\theta_{\p}-\phi_{\p}=2\arcsin(k\sin\pi\lambda_{\p}),
\label{repa5}
\end{equation}
\begin{equation}
\theta_{\p}=\pi\lambda_{\p}+\arcsin(k\sin\pi\lambda_{\p}),
\quad
\phi_{\p}=\pi\lambda_{\p}-\arcsin(k\sin\pi\lambda_{\p}),
\label{repa6}
\end{equation}
which expresses $\theta_{\p}$ and $\phi_{\p}$ in terms of
$\lambda_{\p}$. The other parameter $\gamma_{\p}$ is
defined by
\begin{equation}
\e^{\gamma_{\p}\pm\pi\ii\lambda_{\p}}\equiv
{\e^{\pm\ii\theta_{\p}}+k\over\sqrt{1-k^2}}.
\label{repa7}
\end{equation}
These two expressions are equivalent in view of (\ref{repa2}) and
(\ref{repa4}). Multiplying them and using the second equality in
(\ref{repa3}) we find
\begin{equation}
\e^{2\gamma_{\p}}={1+k^2+2k\cos\theta_{\p}\over1-k^2}=
{\sin\theta_{\p}\over\sin\phi_{\p}}.
\label{repa8}
\end{equation}
Because of (\ref{repa6}) this $\gamma_{\p}$ is also a function
of $\lambda_{\p}$, i.e.
\begin{equation}
\e^{\pm\gamma_{\p}}=
{\sqrt{1-k^2\sin\!\vphantom{k}^2\pi\lambda_{\p}}
\pm k\cos\pi\lambda_{\p}\over\sqrt{1-k^2}}.
\label{repa9}
\end{equation}
From (\ref{intro4}) and (\ref{repa0}) we have
$\mu_{\p}^N=(1+k\,\e^{\theta_{\p}})/k'$ so that
\begin{equation}
\mu_{\p}=
\left({\e^{\ii\theta_{\p}}\sin\theta_{\p}\over
\e^{\ii\phi_{\p}}\sin\phi_{\p}}\right)^{1/2N},
\label{repa10}
\end{equation}
after using (\ref{repa4}) and (\ref{repa7}).

With the help of (\ref{repa0}) and (\ref{repa10}) we can now rewrite
the results (\ref{intro1}) as \cite{Ba}
\begin{eqnarray}
{W_{\pq}(n)\over W_{\pq}(0)}&=&
{\left({\sin\theta_{\p}\sin\phi_{\q}\over
\sin\theta_{\q}\sin\phi_{\p}}\right)}^{n/2N}
\prod^{n}_{j=1}{\sin[\pi(j-\half)/N-(\theta_{\q}-\phi_{\p})/2N]
\over\sin[\pi(j-\half)/N+(\phi_{\q}-\theta_{\p})/2N]},
\nonumber\\&&
\label{repa11}\\
&&\nonumber\\
{\wb_{\pq}(n)\over\wb_{\pq}(0)}&=&
{\left({\sin\theta_{\p}\sin\theta_{\q}\over
\sin\phi_{\p}\sin\phi_{\q}}\right)}^{n/2N}
\prod^{n}_{j=1}
{\sin[\pi(j-1)/N+(\phi_{\q}-\phi_{\p})/2N]
\over\sin[\pi j/N-(\theta_{\q}-\theta_{\p})/2N]}.
\nonumber\\&&
\label{repa12}
\end{eqnarray}
Similarly, using (\ref{repa0}) and (\ref{repa10}), their Fourier
transforms (\ref{intro7}) become\footnote{Eq.\ (7)
and (8) of \cite{AP-inf} have misprints, which can be corrected
by replacing $n$ by $N-n$ in their left-hand sides.}
\begin{eqnarray}
{\wf_{\pq}(n)\over\wf_{\pq}(0)}&=&
{\e^{\ii n(\phi_{\p}-\theta_{\p}+\phi_{\q}-\theta_{\q})/2N}}
\prod^{n}_{j=1}{\sin[\pi(j-1)/N+(\pb_{\q}-\pb_{\p})/2N]
\over\vpOOO\sin[\pi j/N-(\tb_{\q}-\tb_{\p})/2N]},
\nonumber\\&&
\label{repa13}\\
&&\nonumber\\
{\wbf_{\pq}(n)\over\wbf_{\pq}(0)}&=&
{\e^{\ii n(\theta_{\p}-\phi_{\p}-\theta_{\q}+\phi_{\q})/2N}}
\prod^{n}_{j=1}
{\sin[\pi(j-\half)/N-(\pb_{\q}-\tb_{\p})/2N]
\over\vpOOO\sin[\pi(j-\half)/N+(\tb_{\q}-\pb_{\p})/2N]},
\nonumber\\&&
\label{repa14}
\end{eqnarray}
where
\begin{eqnarray}
&\pb_{\p}&=\pi\lambda_{\p}-\ii\gamma_{\p}
=\half(\theta_{\p}+\phi_{\p})-\half
\ii\log{\sin\theta_{\p}\over\sin\phi_{\p}},\nonumber\\
&\tb_{\p}&=\pi\lambda_{\p}+\ii\gamma_{\p}
=\half(\theta_{\p}+\phi_{\p})+\half
\ii\log{\sin\theta_{\p}\over\sin\phi_{\p}}.
\label{repa15}
\end{eqnarray}
By direct substitution we can show that if the weights satisfy
the star-triangle equation (\ref{STE}) then their Fourier transforms
satisfy the star-triangle equation
\begin{eqnarray}
&&{N\over R_{\pqr}}\wbf_{\qr}(a)\,\wf_{\pr}(b)\,\wbf_{\pq}(a+b)
\nonumber\\
&&\qquad=\sum^{N-1}_{d=0}\,\wf_{\pq}(b-d)\,\wbf_{\pr}(a+b-d)\,
\wf_{\qr}(d).
\label{repa16}
\end{eqnarray}
This equation has the exact same form as equation
(\ref{STE}), as can be seen replacing $a\to a-b$, $b\to b-c$,
$a+b\to a-c$, and $c+d\to d$. Therefore, from the proof
\cite{AP-Ta} that the weights (\ref{repa11}), (\ref{repa12})
satisfy (\ref{STE}) we conclude that the weights (\ref{repa13}),
(\ref{repa14}) satisfy (\ref{repa16}).

\par For $\theta_{\p}=\phi_{\p}$, $\theta_{\q}=\phi_{\q}$ we recover
the self-dual Fateev and Zamolodchikov \cite{FZ} solution with
\begin{eqnarray}
&&{\wf(n)\over\wf(0)}={\wb(n)\over\wb(0)}
={\wb(N-n)\over\wb(0)},\nonumber\\
&&\nonumber\\
&&{\wbf(n)\over\wbf(0)}={\wu(n)\over\wu(0)}={\wu(N-n)\over\wu(0)},
\label{repa17}
\end{eqnarray}
which are trigonometric expressions ($g=0$) of the difference
variable $\theta_{\q}-\theta_{\p}$. In this nonchiral special
case the Boltzmann weights depend only on this one parameter,
which is the difference of two rapidity variables. The more
general chiral weights depend on the two rapidity variables
separately, living on a higher-genus curve.

%%%%%%%%%%%%%%%%%%%%%%%%%%%%%%%%%%%%%%%%%%%%%%%%%%%%%%%%%%%%%%%%%%%%%%%%%%%
%  Section 3
%%%%%%%%%%%%%%%%%%%%%%%%%%%%%%%%%%%%%%%%%%%%%%%%%%%%%%%%%%%%%%%%%%%%%%%%%%%
\mypr{}{\setcounter{equation}{0}}
\section{The $N\to\infty$ Limit of the Boltzmann Weights\label{sec-inf}}

In this section we shall obtain the $N\to\infty$ limit of the
Boltzmann weights of the previous section. We shall give explicit
formulae for all three regimes.

%%%%%%%%%%%%%%%%%%%%%%%%%%%%%%%%%%%%%%%%%%%%%%%%%%%%%%%%%%%%%%%%%%%%%%%%%%%
\subsection{General form of the Boltzmann weights\label{susec-bw}}

Note that the Boltzmann weights (\ref{repa11}) and (\ref{repa12}) or
their dual weights (\ref{repa13}) and (\ref{repa14}) all have the product
form
\begin{equation}
{W(n)\over W(0)}=A^{n/N}\prod^{n}_{j=1}
{\sin\left(\vpOO\pi(j+\alpha-1)/N\right)\over
\sin\left(\vpOO\pi(j+\beta-1)/N\right)},
\label{infi0}
\end{equation}
where
\begin{equation}
A={\sin\pi\beta/\sin\pi\alpha},
\label{infi1}
\end{equation}
with $\alpha$ and $\beta$ given constants depending on parameters
$\theta_{\p}$, $\theta_{\q}$, $\phi_{\p}$, and
$\phi_{\q}$ satisfying (\ref{repa1}). Also, the condition on $A$
guarantees that $W(n+N)=W(n)$, using a trivial exercise on complex
exponentials.

More precisely, we have to use in case of (\ref{repa11}) and
(\ref{repa12})
\begin{eqnarray}
\ualpha_{\pq}={1\over2}+{\phi_{\p}-\theta_{\q}\over2\pi},&\quad&
\ubeta_{\pq}={1\over2}+{\phi_{\q}-\theta_{\p}\over2\pi},
\label{infi2}\\
&&\nonumber\\
\balpha_{\pq}={\phi_{\q}-\phi_{\p}\over2\pi},
\phantom{{1\over2}+}&\quad&
\bbeta_{\pq}=1+{\theta_{\p}-\theta_{\q}\over2\pi},
\label{infi3}
\end{eqnarray}
which all four satisfy equations of the form
$\xi_{\pq}+\xi_{\qr}-\xi_{\pr}=\xi_{\qq}$, or
\begin{eqnarray}
&&\ualpha_{\pq}+\ualpha_{\qr}-\ualpha_{\pr}=
\ubeta_{\pq}+\ubeta_{\qr}-\ubeta_{\pr},
\label{infi4}\\
&&\nonumber\\
&&\balpha_{\pq}+\balpha_{\qr}=\balpha_{\pr},\quad
\bbeta_{\pq}+\bbeta_{\qr}=1+\bbeta_{\pr}.
\label{infi5}
\end{eqnarray}
Similarly, we have to use in case of (\ref{repa13}) and (\ref{repa14})
\begin{eqnarray}
\alphaf_{\pq}={\pb_{\q}-\pb_{\p}\over2\pi},
\phantom{{1\over2}+}&\quad&
\betaf_{\pq}=1+{\tb_{\p}-\tb_{\q}\over2\pi},
\label{infi8}\\
&&\nonumber\\
\balphaf_{\pq}=
{1\over2}+{\tb_{\p}-\pb_{\q}\over2\pi},&\quad&
\bbetaf_{\pq}=
{1\over2}+{\tb_{\q}-\pb_{\p}\over2\pi},
\label{infi9}
\end{eqnarray}
satisfying
\begin{eqnarray}
&&\alphaf_{\pq}+\alphaf_{\qr}=\alphaf_{\pr},\quad
\betaf_{\pq}+\betaf_{\qr}=1+\betaf_{\pr},
\label{infi12}\\
&&\nonumber\\
&&\balphaf_{\pq}+\balphaf_{\qr}-\balphaf_{\pr}=
\bbetaf_{\pq}+\bbetaf_{\qr}-\bbetaf_{\pr}.
\label{infi13}
\end{eqnarray}
The four corresponding constants $A$, $\bA$, $\Af$, $\bAf$, as
given in (\ref{infi1}) with the corresponding $\alpha$ and
$\beta$ substituted, are worked out in Appendix \ref{app-A} and
they agree with (\ref{repa11}) through (\ref{repa14}), as was
to be expected.

Important symmetries of weight (\ref{infi0}) are
\begin{eqnarray}
&&{W(n\vert\alpha,\beta)\over W(0\vert\alpha,\beta)}=
{W(0\vert\beta,\alpha)\over W(n\vert\beta,\alpha)}=
{W(n\pm N\vert\alpha,\beta)\over W(0\vert\alpha,\beta)},
\label{infi14a}\\
&&\nonumber\\
&&{W(-n\vert\alpha,\beta)\over W(0\vert\alpha,\beta)}=
{W(N-n\vert\alpha,\beta)\over W(0\vert\alpha,\beta)}=
{W(n\vert1-\beta,1-\alpha)\over W(0\vert1-\beta,1-\alpha)},
\label{infi14}
\end{eqnarray}
which is easily verified from (\ref{infi0}).
This allows us to restrict ourselves to study $W(n)$ only for
$0\le n\le\half N$, while reducing the other case $\half N\le n\le N$
or equivalently $-\half N\le n\le0$ to this case. This symmetry
shows up explicitly in the following, particularly in (\ref{infi18}),
(\ref{infi20}), (\ref{infi22}), and (\ref{B2}).

We can conclude from (\ref{infi14}) that $W(N-n)=W(n)$ for
$\alpha+\beta=1$. Then the chirality disappears and the model
reduces to the model of Fateev and Zamolodchikov \cite{FZ}.

%%%%%%%%%%%%%%%%%%%%%%%%%%%%%%%%%%%%%%%%%%%%%%%%%%%%%%%%%%%%%%%%%%%%%%%%%%%
\subsection{General $N\to\infty$ formula\label{susec-ninf}}

Naively, in the limit $N\to\infty$, we can drop the $\sin$ symbols
in (\ref{infi0}). This leads us to introduce
the function
\begin{eqnarray}
P(n|\alpha,\beta)\equiv
{\Gamma(\alpha+n)\Gamma(\beta)\over\Gamma(\beta+n)\Gamma(\alpha)}
&=\prod^{n}_{j=1}\,{j+\alpha-1\over j+\beta-1}
={(\alpha)_n\over(\beta)_n},\phantom{{}_-}&
\quad\hbox{if }n\ge0,\nonumber\\
&&\nonumber\\
&=\prod^{-n}_{j=1}\,{j-\beta\over j-\alpha}
={(1-\beta)_{-n}\over(1-\alpha)_{-n}},&
\quad\hbox{if }n\le0,
\label{infi15}
\end{eqnarray}
where $\Gamma(x)$ is the Gamma function and
$(x)_n=\Gamma(x+n)/\Gamma(x)$ the Pochhammer symbol \cite{BEMOT,GRJ}.
The finite-$N$ corrections are described by the function
\begin{equation}
S_n(\alpha)\equiv\log\prod_{j=1}^{n}{\sin[\pi(j\!+\!\alpha\!-\!1)/N]\over
\pi(j\!+\!\alpha\!-\!1)/N},
\label{infi16}
\end{equation}
which has an asymptotic expansion derived in Appendix \ref{app-N}. Using
(\ref{A-asym}) there, we immediately have an
asymptotic expansion formula for (\ref{infi0}) in terms of
powers of $1/N$, i.e.\footnote{Formula (12) in \cite{AP-inf}
has the higher orders misprinted and is only correct to the
order needed in the actual $N\to\infty$ limits,
which are presented correctly in \cite{AP-inf}.}
\begin{eqnarray}
\log&&{W(n)\over W(0)}=\log\left[A^{n/N}P(n|\alpha,\beta)\right]
+S_n(\alpha)-S_n(\beta)\nonumber\\&&\nonumber\\
&&=\log\left[A^{n/N}P(n|\alpha,\beta)\right]
+\sum_{l=0}^{\infty}
{\B_{l+1}(\alpha)-\B_{l+1}(\beta)\over(l+1)!}
{\left({\pi\over N}\right)}^l\nonumber\\
&&\qquad\times\left[{\left({\d\over\d z}\right)}^l
{\left.{\log\left({\sin z\over z}\right)}\right\vert}_{z=\pi n/N}-
{\left({\d\over\d z}\right)}^l
{\left.{\log\left({\sin z\over z}\right)}\right\vert}_{z=0}\right],
\label{infi17a}
\end{eqnarray}
where the $\B_m(x)$ are Bernoulli polynomials \cite{BEMOT,GRJ}.
Only the term $l=0$ will be relevant in the limit $N\to\infty$
and the terms $l\ge1$ are finite-$N$ corrections, for which bounds
are derived in Appendix \ref{app-N}. Using $\B_1(x)=x-\half$ and
$\B_2(x)=x^2-x+{1\over6}$ \cite{BEMOT,GRJ} and restricting
ourselves to $l\le1$ we can rewrite (\ref{infi17a}) as
\begin{eqnarray}
{W(n)\over W(0)}&&=A^{n/N}P(n|\alpha,\beta)
\left({\sin(\pi n/N)\over\pi n/N}\right)^{\alpha-\beta}\nonumber\\
&&\nonumber\\
&&\times\exp\left[{\pi(\alpha-\beta)(\alpha+\beta-1)\over2N}
\left(\cot{\pi n\over N}-{N\over\pi n}\right)+\O\big(N^{-2}\big)\right].
\label{infi17}
\end{eqnarray}
The last line of (\ref{infi17}) gives the leading correction
for large $N$ and can be ignored in the limit. We can use
(\ref{infi17}) to study three regimes for the large $N$ limit. We
shall work this out in the following three subsections.

%%%%%%%%%%%%%%%%%%%%%%%%%%%%%%%%%%%%%%%%%%%%%%%%%%%%%%%%%%%%%%%%%%%%%%%%%%%
\subsection{The regime I: $N\to\infty$, $n$ finite\label{susec-w1}}

First we study the limit $N\to\infty$, while
$n$ remains finite. In this case, (\ref{infi17}) results in
\begin{eqnarray}
{W(n)\over W(0)}=&&
{\Gamma(\alpha+n)\Gamma(\beta)\over\Gamma(\beta+n)\Gamma(\alpha)}=
{\Gamma(1-\beta-n)\Gamma(1-\alpha)\over\Gamma(1-\alpha-n)\Gamma(1-\beta)},
\quad-\infty<n<\infty,
\nonumber\\&&
\label{infi18}
\end{eqnarray}
which is just the naive limit $P(n|\alpha,\beta)$ given in (\ref{infi15}).
Using the well-known asymptotic expansion formula of the Gamma function
\cite{BEMOT,GRJ}, we have
\begin{equation}
{\Gamma(z+\alpha)\over\Gamma(z+\beta)}=z^{\alpha-\beta}
\left(1+{(\alpha-\beta)(\alpha+\beta-1)\over2z}+
\O(z^{-2})\right),
\label{infi19}
\end{equation}
for $z\to+\infty$, which is also equation 1.18(4) of \cite{BEMOT}.
Therefore,
\begin{eqnarray}
{W(n)\over W(0)}&=n^{\alpha-\beta}{\Gamma(\beta)\over\Gamma(\alpha)}
\left(1+\O(n^{-1})\right),\phantom{|-||\beta|}&\quad\hbox{for }n\to+\infty,
\nonumber\\
&&\nonumber\\
&=|n|^{\alpha-\beta}{\Gamma(1-\alpha)\over\Gamma(1-\beta)}
\left(1+\O(|n|^{-1})\right),&\quad\hbox{for }n\to-\infty.
\label{infi20}
\end{eqnarray}
This shows that the Boltzmann weights vanish in the limit whenever
$\Re\alpha<\Re\beta$, where $\Re z$ is the real part of $z$.

%%%%%%%%%%%%%%%%%%%%%%%%%%%%%%%%%%%%%%%%%%%%%%%%%%%%%%%%%%%%%%%%%%%%%%%%%%%
\subsection{The regime II: $N,n\to\infty$, $n/N$ finite\label{susec-w2}}

For the second regime we study $N$, $n\to\infty$ such that
\begin{equation}
x\equiv{2\pi n\over N}
\label{infi21}
\end{equation}
remains finite. Consequently, the weights $W(n)$ in (\ref{infi0}),
which originally took $N$ different values and which were periodic
modulo $N$, now depend on the continuous spin values $x$ and they are
periodic modulo $2\pi$.

We can now substitute the asymptotic formula (\ref{infi20}) for
$P(n|\alpha,\beta)$ into (\ref{infi17}), while assuming without
loss of generality $-\half N\le n\le\half N$ and rearranging
the resulting expression. We immediately arrive at
\begin{eqnarray}
W(x)&=W(0)\,{A\vpOO}^{\textstyle{x\over2\pi}}
{\left({N\over\pi}\sin\half x\right)}^{\alpha-\beta}\,
{\Gamma(\beta)\over\Gamma(\alpha)},\phantom{-|\beta|}&
\quad\hbox{if }0<x\le\half\pi,\nonumber\\
&&\nonumber\\
&=W(0)\,{A\vpOO}^{\textstyle{x\over2\pi}}
{\left({N\over\pi}\sin\half|x|\right)}^{\alpha-\beta}\,
{\Gamma(1-\alpha)\over\Gamma(1-\beta)},&
\quad\hbox{if }-\half\pi\le x<0.
\label{infi22}
\end{eqnarray}
This can be summarized as a function periodic modulo $2\pi$, i.e.
\begin{equation}
W(x)=C\,{A\vpOO}^{\textstyle{x\over2\pi}-
\left\lfloor{x\over2\pi}\right\rfloor}\,
{\left\vert\sin{\half x}\right\vert}^{\alpha-\beta},
\label{infi23}
\end{equation}
where $\lfloor x\rfloor$ stands for the largest integer $\le x$ and
\begin{equation}
C=W(0)\,{\left(N\over\pi\right)}^{\alpha-\beta}\,
{\Gamma(\beta)\over\Gamma(\alpha)},\quad
A={\sin\pi\beta\over\sin\pi\alpha}=
{\Gamma(\alpha)\Gamma(1-\alpha)\over\Gamma(\beta)\Gamma(1-\beta)}.
\label{infi24}
\end{equation}
In this regime II, we have to rescale $W(0)$ with a power of $N$
as $N\to\infty$ in order to keep the constant $C$ finite.

For the special case $\alpha+\beta=1$, we have $A=1$ and the
chirality vanishes. This special limit has been mentioned first
by Fateev and Zamolodchikov in \cite{FZ} and is generalized
above to $\alpha+\beta\ne1$.

%%%%%%%%%%%%%%%%%%%%%%%%%%%%%%%%%%%%%%%%%%%%%%%%%%%%%%%%%%%%%%%%%%%%%%%%%%%
\subsection{The regime III: $N,n\to\infty$, $n/N\to0$\label{susec-w3}}

A crossover regime intermediate between regimes I and II appears
when both $N$, $n\to\infty$ such that $n/\varphi(N)\to x$ for some
function $\varphi(N)$ with $\varphi(N)\to\infty$ and
$\varphi(N)/N\to0$. We have
\begin{equation}
W(x)=D\,{A\vpOO}^{-\halfs{\sign}(x)}\,
{\left\vert x\right\vert}^{\alpha-\beta},\quad
-\infty<x<\infty,
\label{infi25}
\end{equation}
which is a chiral generalization of the Boltzmann weight in
Zamolodchikov's Fishnet Model \cite{Za}. Here,
\begin{equation}
D=W(0)\,{\varphi(N)}^{\alpha-\beta}\,
{\Gamma(\beta)\over\Gamma(\alpha)}\,A^{\halfs},
\label{infi26}
\end{equation}
implying again that $W(0)$ need be suitably rescaled in the limit
$N\to\infty$. We note that (\ref{infi25}) is also the asymptotic
large-$n$ behavior (\ref{infi20}) for regime I and the small-$x$
limiting behavior of (\ref{infi23}) in regime II.
The sign function in (\ref{infi25}) arises as coefficients
in (\ref{infi20}) and (\ref{infi22}) differ by a factor $A$ for $n$
positive or negative, see also (\ref{infi24}).

We note that we can reproduce the previously known cases \cite{FZ,Za}
by setting $\alpha+\beta=1$ and $A=1$. Now we have only one condition
(\ref{infi1}) on $A$, i.e.\  $A=\sin\pi\beta/\sin\pi\alpha$. This
provides us with the deformations (\ref{infi18}), (\ref{infi23}), and
(\ref{infi25}), which define integrable field theories with chirality.

%%%%%%%%%%%%%%%%%%%%%%%%%%%%%%%%%%%%%%%%%%%%%%%%%%%%%%%%%%%%%%%%%%%%%%%%%%%
\subsection{Duality of regime I and regime II\label{susec-dual}}

The limiting Boltzmann weights in Regimes I and II are each
other's dual under Fourier duality transformation.
More precisely, if the limiting weights are in Regime II,
their Fourier transforms are in Regime I, and vice versa.
This follows from the way that we have constructed the limits.
However, there is also a direct way to show this, as the infinite
Fourier sum can be performed using a transformation formula of
the Gauss hypergeometric function $\F(a,b;c;x)$ \cite{BEc}.
Thus we obtain a formula for the double-sided hypergeometric
function $\H1$ as defined for example by Slater \cite{Slater},
i.e.\footnote{Start from 2.9 (27) of \cite{BEMOT},
with $a=\alpha$, $b=1$, $c=\beta$, $z=\e^{\ii x}$,
substituting the definitions 2.9 (1), (13), (22). Next use
eqs.\ 2.1.2 (6) and $\F(a,0;c;z)=1$.}
\begin{eqnarray}
&&\hypd{1}{\alpha}{\beta}{\e^{\ii x}}\equiv
\sum_{n=-\infty}^{\infty}{(\alpha)_n\over(\beta)_n}\e^{\ii nx}
=\hypd{1}{1-\beta}{1-\alpha}{\e^{-\ii x}}
\nonumber\\&&
\nonumber\\&&\quad=
\F(\alpha,1;\beta;\e^{\ii x})+\F(1-\beta,1;1-\alpha;\e^{-\ii x})-1
\nonumber\\&&
\nonumber\\&&\quad=
{2^{\beta-\alpha-1}\,\Gamma(1-\alpha)\Gamma(\beta)\over
\Gamma(\beta-\alpha)}\,
\,\e^{\ii(1-\alpha-\beta)(x-\pi)/2}
\,{(\sin\half x)}^{\beta-\alpha-1},
\label{infi27}
\end{eqnarray}
for $0<x<2\pi$ (and periodically extended mod $2\pi$). The inverse
Fourier transform of (\ref{infi27}) corresponds to integral 3.892.1
of \cite{GRJ}, where we need to replace $\nu\mapsto\beta-\alpha$,
$\beta\mapsto\alpha+\beta-1+2n$, $x\mapsto\pi-\half z$,
resulting in
\begin{eqnarray}
{1\over2\pi}\int_0^{2\pi}\d z\,\e^{-\ii nz+\ii(1-\alpha-\beta)z/2}
&&\,{(\sin\half z)}^{\beta-\alpha-1}
={2^{1-\beta+\alpha}\e^{\ii\pi(1-\alpha-\beta)/2}\,(-1)^n
\over(\beta-\alpha)\B(\beta+n,1-\alpha-n)}\nonumber\\
&&={2^{1-\beta+\alpha}\e^{\ii\pi(1-\alpha-\beta)/2}\Gamma(\beta-\alpha)
\,(-1)^n\over\Gamma(\beta+n)\Gamma(1-\alpha-n)},
\label{infi28}
\end{eqnarray}
in agreement with (\ref{infi27}) as
\begin{equation}
(x)_n={\Gamma(x+n)\over\Gamma(x)}=
(-1)^n\,{\Gamma(1-x)\over\Gamma(1-x-n)}={(-1)^n\over(1-x)_{-n}}.
\label{infi28a}
\end{equation}

It is easily verified that, up to a possible overall constant factor,
(\ref{infi27}) and (\ref{infi28}) relate the regime I result
(\ref{infi18}) with the regime II result (\ref{infi23}), provided
an appropriate transformation of the $\alpha$ and $\beta$ parameters
is made. From (\ref{infi2}), (\ref{infi6}), (\ref{repa15}), and
(\ref{infi10}), we find for this duality transformation
\begin{eqnarray}
&&\beta-\alpha-1=\alphaf-\betaf,\nonumber\\
&&\alpha+\beta-1=-{\ii\over\pi}\log \Af,\qquad
\Af={\sin\pi\betaf\over\sin\pi\alphaf},
\label{infi29}
\end{eqnarray}
whereas, from (\ref{infi2}), (\ref{infi6}), (\ref{repa15}), and
(\ref{infi8}), its inverse is given by
\begin{eqnarray}
&&\betaf-\alphaf-1=\alpha-\beta,\nonumber\\
&&\alphaf+\betaf-1={\ii\over\pi}\log A,\qquad
A={\sin\pi\beta\over\sin\pi\alpha}.
\label{infi30}
\end{eqnarray}
The equations (\ref{infi29}) and (\ref{infi30}) differ only by one
minus sign in front of the $\ii$, coming from the corresponding
difference of the finite Fourier transform and its inverse.

%%%%%%%%%%%%%%%%%%%%%%%%%%%%%%%%%%%%%%%%%%%%%%%%%%%%%%%%%%%%%%%%%%%%%%%%%%%
%  Section 4
%%%%%%%%%%%%%%%%%%%%%%%%%%%%%%%%%%%%%%%%%%%%%%%%%%%%%%%%%%%%%%%%%%%%%%%%%%%
\mypr{}{\setcounter{equation}{0}}
\section{The $N\to\infty$ Limit of the Star-Triangle
Equation\label{sec-lim}}

In this section we shall examine the various limits of the
star-triangle equation (\ref{STE}), now we have obtained explicit
prescriptions on how to take the $N\to\infty$ limit
of the Boltzmann weights of the $N$-state chiral Potts model.

%%%%%%%%%%%%%%%%%%%%%%%%%%%%%%%%%%%%%%%%%%%%%%%%%%%%%%%%%%%%%%%%%%%%%%%%%%%
\subsection{Principal domain\label{susec-pd}}

There are several parameter domains that require separate treatment.
But from now on, we shall assume that all $W(n)$ encountered become
vanishingly small whenever $n,N-n\gg 0$. In view of (\ref{infi20})
and (\ref{infi22}), this means that all $\Re(\beta-\alpha)>0$.

From (\ref{infi2}), (\ref{infi3}), (\ref{infi8}), (\ref{infi9}),
(\ref{infi29}), and (\ref{infi30}) we find
\begin{equation}
\beta-\alpha=1-\bbeta+\balpha=
1-\betaf+\alphaf=\bbetaf-\balphaf,
\label{lim1}
\end{equation}
independent of the choice of rapidity variables $\p$, $\q$,
which we have suppressed in (\ref{lim1}). Therefore, we define
our principal domain by the condition
\begin{equation}
0<\Re(\beta-\alpha)<1,
\label{lim2}
\end{equation}
for all occurring ($\alpha$,$\beta$) pairs. This condition is
easily satisfied, even with $\ualpha_{\pq}$, $\ubeta_{\pq}$,
$\ualpha_{\pr}$, $\ubeta_{\pr}$, $\ualpha_{\qr}$, $\ubeta_{\qr}$,
and their barred versions, all being real.

The summation over $d$ in (\ref{STE}) has to be split in several
pieces as we must choose to which of the three $N\to\infty$
regimes each of the three weights in the summand belongs. We shall
see that under condition (\ref{lim2}) all pieces but one can be ignored
and that the three types of large-$N$ behavior I, II, or III do not mix:
If we take the three spin states $a$, $b$, and $c$ in (\ref{STE}) such
that all three weights in the right-hand side of (\ref{STE}) are in the
same regime, the dominant part of the sum over spin state $d$ comes
from the piece with all three weights in the left-hand side of
(\ref{STE}) being in the identical regime.

%%%%%%%%%%%%%%%%%%%%%%%%%%%%%%%%%%%%%%%%%%%%%%%%%%%%%%%%%%%%%%%%%%%%%%%%%%%
\subsection{The constant $R_{\pqr}$\label{susec-R}}

The next calculation to be done is the evaluation of the large-$N$
limit of the constants $R_{\pqr}$ or $F_{\pq}$ given in (\ref{STER}).
We break this up in several steps.

First, we can rewrite the $F_{\pq}$ in (\ref{STER}) as
\begin{equation}
F_{\pq}=N{\wbf_{\pq}(0)\over W_{\pq}(0)}
\exp({\bLf_{\pq}-L_{\pq}}),
\label{lim3}
\end{equation}
where
\begin{equation}
L_{\pq}\equiv{1\over N}\sum_{l=1}^{N}\log{W_{\pq}(l)\over W_{\pq}(0)},
\qquad\bLf_{\pq}\equiv
{1\over N}\sum_{l=1}^{N}\log{\wbf_{\pq}(l)\over\wbf_{\pq}(0)}.
\label{lim4}
\end{equation}
Here both the $L_{\pq}$ and the $\bLf_{\pq}$ can be evaluated
in an identical fashion under condition (\ref{lim2}), substituting
the regime II asymptotic form (\ref{infi23}) and (\ref{infi24}),
while replacing the sum by an integral. Therefore, in an obvious
simplification of notation suppressing the rapidity subscripts,
$L_{\pq}$ becomes
\begin{eqnarray}
L&&=\log\left[{(N/\pi)}^{\alpha-\beta}\,
\Gamma(\beta)/\Gamma(\alpha)\right]\nonumber\\
&&\qquad+{1\over N}\sum_{l=1}^{N}\left[
\Big((n/N)-\lfloor n/N\rfloor\Big)\log A+
(\alpha-\beta)\log\sin|\pi n/N|\,\right]\nonumber\\
&&\nonumber\\
&&\approx\log\left[{(N/\pi)}^{\alpha-\beta}\,
\Gamma(\beta)/\Gamma(\alpha)\right]\nonumber\\
&&\qquad+{1\over2\pi}\int_{0}^{2\pi}\d x\,\left[
\Big((x/2\pi)-\lfloor x/2\pi\rfloor\Big)\log A+
(\alpha-\beta)\log\sin|\half x|\,\right]\nonumber\\
&&\nonumber\\
&&=\log\left[{(N/2\pi)}^{\alpha-\beta}\,A^{1/2}\,
\Gamma(\beta)/\Gamma(\alpha)\right],
\qquad \hbox{for }N\to\infty,
\label{lim5}
\end{eqnarray}
where the elementary integral 4.224.3 of \cite{GRJ} has been used.
Similarly,
\begin{equation}
\bLf\approx
\log\Bigl[{(N/2\pi)}^{\balphaf-\bbetaf}\,\mbox{$\bAf$}^{1/2}\,
\Gamma(\bbetaf)/\Gamma(\balphaf)\Bigr],
\quad \hbox{for }N\to\infty.
\label{lim6}
\end{equation}
It is easily checked that the corrections to (\ref{lim5}) and
 (\ref{lim6}) are irrelevant in the large-$N$ limit.

Similarly, $\wbf(0)/\wb(0)$ is also dominated by the regime II
contribution given by (\ref{infi23}) and (\ref{infi24}), i.e.
\begin{eqnarray}
{\wbf(0)\over\wb(0)}&&={1\over N}\sum_{l=1}^{N}
\,\Bigl[{(N/\pi)}^{\balpha-\bbeta}\,
\Gamma(\bbeta)/\Gamma(\balpha)\Bigr]\,
\bA^{(n/N)-\lfloor n/N\rfloor}
{|\sin{\pi n/N}|}^{\balpha-\bbeta}\nonumber\\
&&\nonumber\\
&&\approx\Bigl[{(N/\pi)}^{\balpha-\bbeta}\,
\Gamma(\bbeta)/\Gamma(\balpha)\Bigr]\,
{1\over2\pi}\int_{0}^{2\pi}\d x\,\bA^{x/2\pi}
{(\sin\half x)}^{\balpha-\bbeta}\nonumber\\
&&\nonumber\\
&&={(N/2\pi)}^{\balpha-\bbeta}\,\bA^{1/2}\,
{\Gamma(\bbeta)\Gamma(1-\bbeta+\balpha)\over
\Gamma(\balpha)\Gamma(\bbetaf)
\Gamma(1-\balphaf)},
\quad \hbox{for }N\to\infty.
\label{lim7}
\end{eqnarray}
We have used the integral (\ref{infi28}) for $n=0$,
after substituting the barred version of (\ref{infi30}), i.e.
\begin{equation}
\balpha-\bbeta=\bbetaf-\balphaf-1,
\quad\log\bA=\pi\ii\,(1-\balphaf-\bbetaf).
\label{lim8}
\end{equation}
Then the last line of (\ref{lim7}) is obtained using
(\ref{lim8}) once more. As before in (\ref{lim5}) and (\ref{lim6}),
(\ref{lim7}) gives the coefficient of the leading $N$-power as
$N\to\infty$. Correction terms can be obtained, but they will
not be needed.

We can now substitute (\ref{lim4})--(\ref{lim8}) into (\ref{lim3})
and simplify the result. This is worked out in Appendix \ref{app-F}
and the result is
\begin{equation}\phantom{_pq}
F_{\pq}=
\Bigl[{({N/2\pi})}_{\vphantom{\pq}}^{\strut\ubeta_{\pq}-\ualpha_{\pq}}
\,A_{\pq}^{\strut\,-1/2}\Bigr]\,{\wb_{\pq}(0)\over W_{\pq}(0)}\,
{\Gamma(\ualpha_{\pq})\Gamma(\bbeta_{\pq})\Gamma(1-\balpha_{\pq})
\over\Gamma(\ubeta_{\pq})\Gamma(\bbeta_{\pq}-\balpha_{\pq})},
\label{lim9}
\end{equation}
giving us the desired expression for
$R_{\pqr}={F_{\pq}F_{\qr}/F_{\pr}}$. We note that the factor in
square brackets in (\ref{lim9})  cancels out in view of (\ref{infi4})
and (\ref{infi6}). Hence,
\begin{equation}
\lim_{N\to\infty}R_{\pqr}={\wb_{\pq}(0)W_{\pr}(0)\wb_{\qr}(0)\over
W_{\pq}(0)\wb_{\pr}(0)W_{\qr}(0)}\,r_{\pqr}^{\infty},
\label{lim10}
\end{equation}
where
\begin{equation}
r_{\pqr}^{\infty}\equiv{f_{\pq}f_{\qr}\over f_{\pr}},\quad
f_{\pq}\equiv
{\Gamma(\ualpha_{\pq})\Gamma(\bbeta_{\pq})\Gamma(1-\balpha_{\pq})
\over\Gamma(\ubeta_{\pq})\Gamma(\bbeta_{\pq}-\balpha_{\pq})}.
\label{lim11}
\end{equation}
We see that $R_{\pqr}$ is independent of $N$ in leading order
in the large-$N$ limit. Finite-$N$ corrections can be worked out
but are not needed here, as we shall only consider the actual
$N\to\infty$ limit.

Note that $f_{\pq}$ and $r_{\pqr}^{\infty}$ are invariant under
$\balpha_{\pq}\mapsto1-\bbeta_{\pq}$,
$\bbeta_{\pq}\mapsto1-\balpha_{\pq}$,
whereas $F_{\pq}$ and $r_{\pqr}^{\infty}$
are invariant under $\ualpha_{\pq}\mapsto1-\ubeta_{\pq}$,
$\ubeta_{\pq}\mapsto1-\ualpha_{\pq}$, but $f_{\pq}$ is not.

%%%%%%%%%%%%%%%%%%%%%%%%%%%%%%%%%%%%%%%%%%%%%%%%%%%%%%%%%%%%%%%%%%%%%%%%%%%
\subsection{Regime I\label{susec-ste1}}

We can first consider the large-$N$ limit of (\ref{STE}), while
keeping $a-b$, $a-c$, and thus also $b-c$ finite. Without loss of
generality, we can then restrict ourselves to considering the case
$a$, $b$, and $c$ finite. The three Boltzmann weights in the
right-hand side of (\ref{STE}) take the regime-I form (\ref{infi18}).
The sum over $d$ in (\ref{STE}) needs to be split up in the limit.
In one part $|d|$ remains finite but can become arbitrarily large
and the three weights in the left-hand side of (\ref{STE}) also
take the regime-I form. The summand of (\ref{STE}) then decays
as ${|d|}^{\kappa}$ for $|d|\to\infty$. From (\ref{infi18}), and
using (\ref{infi5}) and (\ref{lim1}), we see that
\begin{eqnarray}
\kappa&&=\balpha_{\qr}-\bbeta_{\qr}+\ualpha_{\pr}-\ubeta_{\pr}
+\balpha_{\pq}-\bbeta_{\pq}\nonumber\\
&&=\balpha_{\qr}-\bbeta_{\qr}+\bbeta_{\pr}-\balpha_{\pr}-1
+\balpha_{\pq}-\bbeta_{\pq}=-2,
\label{lim12}
\end{eqnarray}
so that the sum over $d$ converges as $\sum{|d|}^{-2}$.

In the part of the sum for which $|d|/N$ does not tend to zero,
the three weights in the left-hand side of (\ref{STE}) belong to
regime II and now we can use (\ref{infi22}) to show that the
summand scales as $N^{\kappa}=N^{-2}$. As the total sum has $N$
terms, this contribution vanishes in the limit.

The contribution of the crossover regime III connecting
regimes I and II also vanishes. To show this in more detail,
we can split the sum over $d$ in a piece $-N^{2/3}<d\le N^{2/3}$,
and a piece $N^{2/3}<d\le N-N^{2/3}$. In the first piece,
the summand is bounded by ${|d|}^{-2}$-behavior as shown
above and in (\ref{infi25}) and (\ref{infi26}); therefore,
the error made by replacing the sum with $-N^{2/3}<d\le N^{2/3}$
by a sum $-\infty<d\le\infty$ vanishes as $\O(N^{-2/3})$.
For the second piece we can use (\ref{infi22}), (\ref{infi25})
and (\ref{infi26}) to show that the summand is $\O(N^{2\kappa/3})$
and with less than $N$ terms of order $\O(N^{-4/3})$
its contribution vanishes as $\O(N^{-1/3})$.

To summarize, all six Boltzmann weights must take the
regime-I form (\ref{infi18})
\begin{equation}
W_{\pq}(n)={{(\ualpha_{\pq})}_n\over{(\ubeta_{\pq})}_n},
\qquad\wb_{\pq}(n)={{(\balpha_{\pq})}_n\over{(\bbeta_{\pq})}_n},
\label{lim13}
\end{equation}
solving the star-triangle equation
\begin{eqnarray}
&&\sum^{\infty}_{d=-\infty}\,
\wb_{\qr}(b-d)\,W_{\pr}(a-d)\,\wb_{\pq}(d-c)
\nonumber\\
&&\qquad=r_{\pqr}^{\infty}\,
W_{\pq}(a-b)\,\wb_{\pr}(b-c)\,W_{\qr}(a-c),
\label{ste0}
\end{eqnarray}
provided the $\alpha$ and $\beta$ parameters satisfy (\ref{infi2}),
(\ref{infi3}), and (\ref{repa2}). In (\ref{lim13}) and (\ref{ste0})
we have chosen the normalization $W_{\pq}(0)=\wb_{\pq}(0)=1$.

Equation (\ref{ste0}) is related to the Dougall--Ramanujan
identity \cite{BEd}. We shall return to this in the next section.

%%%%%%%%%%%%%%%%%%%%%%%%%%%%%%%%%%%%%%%%%%%%%%%%%%%%%%%%%%%%%%%%%%%%%%%%%%%
\subsection{Regime II\label{susec-ste2}}

We can next consider the large-$N$ limit of (\ref{STE}), while
keeping $(a-b)/N$, $(a-c)/N$, and $(b-c)/N$ fixed and nonzero.
The three Boltzmann weights in the right-hand side of (\ref{STE})
now take the regime-II form (\ref{infi22}). Again, the sum over
$d$ in (\ref{STE}) needs to be split up in this limit. Now the
dominant part is the one with all three weights in the left-hand
side of (\ref{STE}) belonging to regime II. Only when $d$ is close
to $a$, $b$, or $c$, one of the weights can be of the form of
regime I or III. Because of the principal domain condition
(\ref{lim2}) these contributions can be ignored in the large-$N$
limit and the sum can be replaced by an integral as is done twice
in subsection \ref{susec-R}.

In the previous section we have seen that the summand in the
left-hand side of (\ref{STE}) scales as $N^{-2}$ in the large-$N$
limit, when the three weights belong to regime II. The right-hand
side of (\ref{STE}) now scales as $N^{\bkappa}=N^{-1}$, since
\begin{eqnarray}
\bkappa&&=\ualpha_{\pq}-\ubeta_{\pq}+\balpha_{\pr}-\bbeta_{\pr}
+\ualpha_{\qr}-\ubeta_{\qr}\nonumber\\
&&=\ualpha_{\pq}-\ubeta_{\pq}+\ubeta_{\pr}-\ualpha_{\pr}-1
+\ualpha_{\qr}-\ubeta_{\qr}=-1.
\label{lim14}
\end{eqnarray}
Therefore, it is natural to multiply (\ref{STE}) by $N$ and to
replace $N^{-1}\sum_d$ by $(2\pi)^{-1}\int\d x$, so that the
star-triangle equation becomes
\begin{eqnarray}
{1\over2\pi}
&&\int^{2\pi}_0\d w\,
\wb_{\qr}(y-w)\,W_{\pr}(x-w)\,\wb_{\pq}(w-z)\nonumber\\
&&\qquad=R_{\pqr}^{\infty}\,
W_{\pq}(x-y)\,\wb_{\pr}(y-z)\,W_{\qr}(x-z),
\label{ste1}
\end{eqnarray}
as also follows after suitable rescalings of the Boltzmann weights
and $R_{\pqr}$, i.e.
\begin{equation}
R_{\pqr}^{\infty}\equiv\lim_{N\to\infty}N^{-1}R_{\pqr}
\label{lim15}
\end{equation}
and the two equations (\ref{lim17}) below. Equation (\ref{ste1})
has the solution
\begin{eqnarray}
W_{\pq}(x)&=&{\e\vpOO}^{(\gamma_{\p}-\gamma_{\q})
({\textstyle{x\over2\pi}-\left\lfloor{x\over2\pi}\right\rfloor})}\,
{\left\vert\sin{\half x}\right\vert}^{\lambda_{\p}-\lambda_{\q}},
\nonumber\\
\wb_{\pq}(x)&=&{\e\vpOO}^{(\gamma_{\p}+\gamma_{\q})
({\textstyle{x\over2\pi}-\left\lfloor{x\over2\pi}\right\rfloor})}\,
{\left\vert\sin{\half x}\right\vert}^{\lambda_{\q}-\lambda_{\p}-1},
\label{ste2}
\end{eqnarray}
using the normalization $C_{\pq}=\bC_{\pq}=1$ in (\ref{infi23}),
together with (\ref{infi2}), (\ref{infi3}), (\ref{repa4}),
(\ref{infi6}), and (\ref{infi7}), which imply
\begin{eqnarray}
&&A_{\pq}=\e^{\gamma_{\p}-\gamma_{\q}},\qquad
\ualpha_{\pq}-\ubeta_{\pq}=\lambda_{\p}-\lambda_{\q},\nonumber\\
&&\bA_{\pq}=\e^{\gamma_{\p}+\gamma_{\q}},\qquad
\balpha_{\pq}-\bbeta_{\pq}=\lambda_{\q}-\lambda_{\p}-1.
\label{lim16}
\end{eqnarray}
Here $\gamma_{\p}$ and $\lambda_{\p}$ are related by (\ref{repa9}).
If $\lambda_{\p}<\lambda_{\q}<\lambda_{\r}<1+\lambda_{\p}$ all six
Boltzmann weights in (\ref{ste1}) are real and positive and the
parameters are in the principal domain.

According to (\ref{infi24}), the condition $C_{\pq}=\bC_{\pq}=1$
implies
\begin{equation}
W_{\pq}(0)=
{\left(N\over\pi\right)}^{\ubeta_{\pq}-\ualpha_{\pq}}\,
{\Gamma(\ualpha_{\pq})\over\Gamma(\ubeta_{\pq})},\quad
\wb_{\pq}(0)=
{\left(N\over\pi\right)}^{\bbeta_{\pq}-\balpha_{\pq}}\,
{\Gamma(\balpha_{\pq})\over\Gamma(\bbeta_{\pq})}.
\label{lim17}
\end{equation}
Substituting this in (\ref{lim15}), while using (\ref{lim10}) and
(\ref{lim11}), we arrive at
\begin{equation}
R_{\pqr}^{\infty}=
{\tilde f_{\pq}\tilde f_{\qr}\over\tilde f_{\pr}},\quad
\tilde f_{\pq}\equiv{1\over\pi}\,
{\Gamma(\ubeta_{\pq})\Gamma(\balpha_{\pq})\over
\Gamma(\ualpha_{\pq})\Gamma(\bbeta_{\pq})}\,f_{\pq}=
{\Gamma(\balpha_{\pq})\Gamma(1-\balpha_{\pq})
\over\pi\Gamma(\bbeta_{\pq}-\balpha_{\pq})},
\label{ste3}
\end{equation}
where the $\pi$ and $N$ factors have been redistributed
using (\ref{infi4}) and (\ref{infi5}).

It can be shown that the Boltzmann weights obtained by
dropping the integral part $\lfloor x/2\pi\rfloor$ in
(\ref{ste2}) also satisfy the same star-triangle equation
(\ref{ste1}). The resulting chiral solution can be viewed as the
nonchiral Fateev--Zamolodchikov large-$N$ solution \cite{FZ}
with a site-dependent gauge transformation.

Since $W_{\pq}(x)$ and $\wb_{\pq}(x)$ as given in (\ref{ste2}) are
now functions of $x$ periodic modulo $2\pi$, their Fourier transforms
\begin{eqnarray}
&&\wf_{\pq}(j)=
{1\over2\pi}\int_{0}^{2\pi}\d x\,\e^{-\ii jx}\,W_{\pq}(x),
\quad
\wbf_{\pq}(j)=
{1\over2\pi}\int_{0}^{2\pi}\d x\,\e^{-\ii jx}\,\wb_{\pq}(x),
\nonumber\\&&
\label{ste4}
\end{eqnarray}
are over all integer value $j$, ranging from $-\infty$ to $\infty$.
Substituting (\ref{ste4}) into (\ref{ste1}), we find that these
Fourier transforms satisfy the star-triangle equation
\begin{eqnarray}
&&{1\over R_{\pqr}^{\infty}}\,
\wbf_{\qr}(a-b)\,\wf_{\pr}(b-c)\,\wbf_{\pq}(a-c)\nonumber\\
&&\quad=\sum^{\infty}_{d=-\infty}\,
\wf_{\pq}(b-d)\,\wbf_{\pr}(a-d)\,\wf_{\qr}(d-c),
\label{ste5}
\end{eqnarray}
in which the sum is over all integer values of $d$.

%%%%%%%%%%%%%%%%%%%%%%%%%%%%%%%%%%%%%%%%%%%%%%%%%%%%%%%%%%%%%%%%%%%%%%%%%%%
\subsection{Regime III\label{susec-ste3}}

As the final case, we can consider the large-$N$ limit of
(\ref{STE}), while keeping $(a-b)/\varphi(N)$, $(a-c)/\varphi(N)$,
and $(b-c)/\varphi(N)$ fixed at finite and nonzero values.
The three Boltzmann weights in the right-hand side of (\ref{STE})
now take the regime-III form (\ref{infi25}). Also in this case,
the sum over $d$ in (\ref{STE}) needs to be split up. Now the
dominant part is the one with all three weights in the left-hand
side of (\ref{STE}) belonging to regime III. Only when $d$ is close
to $a$, $b$, or $c$, one of these weights can be of the form of
regime I, while the other two weights are of the form of regime III.
Because of the principal domain condition (\ref{lim2}) these
regime-I contributions can be ignored in the large-$N$ limit and
the sum can be replaced by an integral as is done twice already in
subsections \ref{susec-R} and \ref{susec-ste2}. Also, the three weights
in the summand of (\ref{STE}) can only simultaneously take the form
(\ref{infi22}) of regime II. This contribution is $\O(N^{-1})$ as
in subsection \ref{susec-ste1} and can be ignored.

Hence, we can prove a star-triangle equation of the form
(\ref{ste1}), but with integration over $(-\infty,+\infty)$, i.e.
\begin{eqnarray}
&&\int^{+\infty}_{-\infty}\d w\,
\wb_{\qr}(y-w)\,W_{\pr}(x-w)\,\wb_{\pq}(w-z)\nonumber\\
&&\qquad=\hat R_{\pqr}^{\infty}\,
W_{\pq}(x-y)\,\wb_{\pr}(y-z)\,W_{\qr}(x-z),
\label{ste6}
\end{eqnarray}
where
\begin{equation}
\hat R_{\pqr}^{\infty}\equiv\lim_{N\to\infty}\varphi(N)^{-1}R_{\pqr}.
\label{lim18}
\end{equation}
Equation (\ref{ste6}) has the solution
\begin{eqnarray}
W_{\pq}(x)&=&{\e\vpOO}^{-\halfs(\gamma_{\p}-\gamma_{\q})\sign(x)}\,
{\vert x\vert}^{\lambda_{\p}-\lambda_{\q}},
\nonumber\\
&&\nonumber\\
\wb_{\pq}(x)&=&{\e\vpOO}^{-\halfs(\gamma_{\p}+\gamma_{\q})\sign(x)}\,
{\vert x\vert}^{\lambda_{\q}-\lambda_{\p}-1},
\label{ste7}
\end{eqnarray}
using the normalization $D_{\pq}=\bD_{\pq}=1$ in (\ref{infi25}),
together with (\ref{lim16}).
Again, $\gamma_{\p}$ and $\lambda_{\p}$ are related by (\ref{repa9}).
If $\lambda_{\p}<\lambda_{\q}<\lambda_{\r}<1+\lambda_{\p}$ all six
Boltzmann weights in (\ref{ste6}) are real and positive.

According to (\ref{infi26}), the condition $D_{\pq}=\bD_{\pq}=1$
implies
\begin{equation}
W_{\pq}(0)=
{\varphi(N)}^{\ubeta_{\pq}-\ualpha_{\pq}}\,
{\Gamma(\ualpha_{\pq})\over\Gamma(\ubeta_{\pq})},\quad
\wb_{\pq}(0)=
{\varphi(N)}^{\bbeta_{\pq}-\balpha_{\pq}}\,
{\Gamma(\balpha_{\pq})\over\Gamma(\bbeta_{\pq})}.
\label{lim20}
\end{equation}
Substituting this in (\ref{lim18}), while using (\ref{lim10}) and
(\ref{lim11}), we arrive at
\begin{equation}
\hat R_{\pqr}^{\infty}=
{\hat f_{\pq}\hat f_{\qr}\over\hat f_{\pr}},\quad
\hat f_{\pq}\equiv
{\Gamma(\ubeta_{\pq})\Gamma(\balpha_{\pq})\over
\Gamma(\ualpha_{\pq})\Gamma(\bbeta_{\pq})}\,f_{\pq}=
{\Gamma(\balpha_{\pq})\Gamma(1-\balpha_{\pq})
\over\Gamma(\bbeta_{\pq}-\balpha_{\pq})},
\label{ste8}
\end{equation}
where the $\varphi(N)$ factors have been redistributed
using (\ref{infi4}) and (\ref{infi5}).

It can be shown that the nonchiral Boltzmann weights
obtained by setting $\gamma_{\p}=\gamma_{\q}=0$ in (\ref{ste7})
satisfy the same star-triangle equation (\ref{ste6}). The
resulting solution can be viewed as the Fishnet Model
of Zamolodchikov \cite{Za}. It also generalizes Symanzik's
conformal integral \cite{Symanzik}, which has been used by
Zamolodchikov to prove the star-triangle equation for the
Fishnet Model and which has also provided the proof for the
$N=\infty$ Fateev-Zamolodchikov model \cite{FZ} via a conformal
transformation, $\xi=\tan\half w$.

%%%%%%%%%%%%%%%%%%%%%%%%%%%%%%%%%%%%%%%%%%%%%%%%%%%%%%%%%%%%%%%%%%%%%%%%%%%
\subsection{Remark on $\R$-matrix\label{susec-rmat}}

Following our joint work with Baxter \cite{BPA}, we can make an
$\R$-matrix by taking the product of four weights in any of the
above regimes I, II, or III. Taking four weights of type I, we have
\begin{equation}
\R(a,b,c,d)=\wb_{\p_1\q_1}(a-c)\,W_{\p_1\q_2}(c-b)\,
\wb_{\p_2\q_2}(d-b)\,W_{\p_2\q_1}(a-d).
\label{ste9}
\end{equation}
Similarly, we can also take four weights of type II or III.
Then any such infinite-dimensional $\R$-matrix satisfies the usual
Yang-Baxter equation. But these solutions are very different from
those of \cite{Ba-SOS,FZ-rot,Gaudin,SU,Shibukawa}.

%%%%%%%%%%%%%%%%%%%%%%%%%%%%%%%%%%%%%%%%%%%%%%%%%%%%%%%%%%%%%%%%%%%%%%%%%%%
%  Section 5
%%%%%%%%%%%%%%%%%%%%%%%%%%%%%%%%%%%%%%%%%%%%%%%%%%%%%%%%%%%%%%%%%%%%%%%%%%%
\mypr{}{\setcounter{equation}{0}}
\section{Two-Sided Hypergeometric Sum\label{sec-hyp}}

In this section we shall rewrite the star-triangle equation
(\ref{ste0}) with solution (\ref{lim13}) as a new double-sided
hypergeometric identity.

%%%%%%%%%%%%%%%%%%%%%%%%%%%%%%%%%%%%%%%%%%%%%%%%%%%%%%%%%%%%%%%%%%%%%%%%%%%
\subsection{More symmetric star-triangle equation\label{susec-symste}}

The star-triangle equation (\ref{ste0}) can be written in a
more symmetric form with permutation symmetry among the Boltzmann
weights in each of the two sides.

We start by applying (\ref{infi14}) to the $\wb_{\pq}$ in
(\ref{ste0}), followed by substituting $d=-n$ and and replacing
all six weights by (\ref{lim13}). We arrive at
\begin{eqnarray}
&&\sum^{\infty}_{n=-\infty}\,
{{(\balpha_{\qr})}_{b+n}\,{(\ualpha_{\pr})}_{a+n}\,
{(1-\bbeta_{\pq})}_{c+n}\over
{(\bbeta_{\qr})}_{b+n}\,{(\ubeta_{\pr})}_{a+n}\,
{(1-\balpha_{\pq})}_{c+n}}
=r_{\pqr}^{\infty}\,
{{(\ualpha_{\pq})}_{a-b}\,{(\balpha_{\pr})}_{b-c}\,
{(\ualpha_{\qr})}_{a-c}\over
{(\ubeta_{\pq})}_{a-b}\,{(\bbeta_{\pr})}_{b-c}\,
{(\ubeta_{\qr})}_{a-c}},
\nonumber\\&&
\label{hyp1}
\end{eqnarray}
This result seems to invite us to introduce a more
symmetric notation. For the three external spin
states in the star-triangle equation we write
\begin{equation}
m_1\equiv a,\qquad m_2\equiv b,\qquad m_3\equiv c.
\label{hyp2}
\end{equation}
It is logical to associate $\pr$ with 1, $\qr$ with 2, and
$\pq$ with 3. The other quantities in (\ref{hyp1}) are then
rewritten as\footnote{The parameters $x_j$ and $y_j$,
for $j=1,2,3$, should not be confused with the rapidity
variables $x_{\p}$ and $y_{\p}$ of (\ref{intro1}).}
\begin{eqnarray}
x_1\equiv\ualpha_{\pr},\qquad
&x_2\equiv\balpha_{\qr},\qquad
&x_3\equiv1-\bbeta_{\pq},\nonumber\\
y_1\equiv\ubeta_{\pr},\qquad
&y_2\equiv\bbeta_{\qr},\qquad
&y_3\equiv1-\balpha_{\pq},\nonumber\\
u_1\equiv\balpha_{\pr},\qquad
&u_2\equiv\ualpha_{\qr},\qquad
&u_3\equiv\ualpha_{\pq},\nonumber\\
v_1\equiv\bbeta_{\pr},\qquad
&v_2\equiv\ubeta_{\qr},\qquad
&v_3\equiv\ubeta_{\pq},
\label{hyp3}\\
&&\nonumber\\
f_1\equiv f_{\pr},\qquad
&f_2\equiv f_{\qr},\qquad
&f_3\equiv f_{\pq},
\label{hyp4}
\end{eqnarray}
where the $f_j$ for $j=1,2,3$ are given by (\ref{lim11}), i.e.
\begin{eqnarray}
&&f_1={\Gamma(x_1)\Gamma(v_1)\Gamma(1-u_1)
\over\Gamma(y_1)\Gamma(v_1-u_1)},\nonumber\\
&&f_2={\Gamma(u_2)\Gamma(y_2)\Gamma(1-x_2)
\over\Gamma(v_2)\Gamma(y_2-x_2)},\quad
f_3={\Gamma(u_3)\Gamma(y_3)\Gamma(1-x_3)
\over\Gamma(v_3)\Gamma(y_3-x_3)},
\label{hyp5}
\end{eqnarray}
and $r_{\pqr}^{\infty}=f_2f_3/f_1$.
Therefore, (\ref{hyp1}) takes the much more symmetric form
\begin{eqnarray}
&&\sum^{\infty}_{n=-\infty}\,
{{(x_1)}_{m_1+n}\,{(x_2)}_{m_2+n}\,{(x_3)}_{m_3+n}\over
{(y_1)}_{m_1+n}\,{(y_2)}_{m_2+n}\,{(y_3)}_{m_3+n}}\nonumber\\
&&\nonumber\\
&&\qquad={f_2f_3\over f_1}\,
{{(u_1)}_{m_1-m_2}\,{(u_2)}_{m_2-m_3}\,{(u_3)}_{m_1-m_3}\over
{(v_1)}_{m_1-m_2}\,{(v_2)}_{m_2-m_3}\,{(v_3)}_{m_1-m_3}}.
\label{hyp6}
\end{eqnarray}
The quantities (\ref{hyp3}) are not independent. From
(\ref{infi2}) and (\ref{infi3}) we find that the six
variables $u_j$ and $v_j$ depend in a linear and symmetric
fashion on the six variables $x_j$ and $y_j$, i.e.
\begin{eqnarray}
&u_1=1+x_2-y_3,\qquad&v_1=y_2-x_3,\nonumber\\
&u_2=1+x_1-y_3,\qquad&v_2=y_1-x_3,\nonumber\\
&u_3=1+x_1-y_2,\qquad&v_3=y_1-x_2,
\label{hyp7}
\end{eqnarray}
so that (\ref{hyp6}) simplifies further to
\begin{equation}
\sum^{\infty}_{n=-\infty}\,\prod_{j=1}^3
{{(x_j)}_{m_j+n}\over{(y_j)}_{m_j+n}}
={f_2f_3\over f_1}\,\myunderset{1\le i<j\le3}{\prod\prod}
{{(1+x_i-y_j)}_{m_i-m_j}\over{(y_i-x_j)}_{m_i-m_j}}.
\label{hyp8}
\end{equation}
Since (\ref{infi28a}) implies
\begin{equation}
{{(1+x_i-y_j)}_{m_i-m_j}\over{(y_i-x_j)}_{m_i-m_j}}=
{{(1+x_j-y_i)}_{m_j-m_i}\over{(y_j-x_i)}_{m_j-m_i}},
\label{hyp9}
\end{equation}
the symmetry in (\ref{hyp8}) is even larger than is manifested there,
namely the full permutation symmetry group $\hbox{S}_3$.

%%%%%%%%%%%%%%%%%%%%%%%%%%%%%%%%%%%%%%%%%%%%%%%%%%%%%%%%%%%%%%%%%%%%%%%%%%%
\subsection{The two conditions to be satisfied\label{susec-cond}}

The integrable chiral Potts model solution (\ref{intro1}) of the
star-triangle equation (\ref{STE}) has not six free parameters,
but only four independent ones, namely the rapidity variables
$x_{\p}$, $x_{\q}$, $x_{\r}$, and the modulus parameter $k$.
Therefore, two conditions must be imposed on the new variables
$x_j$ and $y_j$ with $j=1,2,3$.

The first relation is a linear relation, which is a direct
consequence of (\ref{lim12}) and (\ref{hyp3}). It reads
\begin{equation}
x_1+x_2+x_3+2=y_1+y_2+y_3.
\label{hyp10}
\end{equation}
This is in fact the Saalsch\"utz condition, which plays such
an important role in the theory of hypergeometric functions
\cite{BEMOT,Slater}.

The second relation can be found from (\ref{lim16}), which
is derived in Appendix \ref{app-A} and holds for all $N$.
From (\ref{lim16}) we can derive
\begin{equation}
A_{\pr}\bA_{\qr}=\bA_{\pq}=\bA_{\qp}=A_{\qr}\bA_{\pr},\quad
A_{\pq}A_{\qr}=A_{\pr},\quad A_{\pq}A_{\qp}=1.
\label{hyp11}
\end{equation}
From the first equality in (\ref{hyp11}) and (\ref{hyp3}), we
then find
\begin{equation}
\sin\pi x_1\,\sin\pi x_2\,\sin\pi x_3=
\sin\pi y_1\,\sin\pi y_2\,\sin\pi y_3,
\label{hyp12}
\end{equation}
and we may call this the ``periodicity condition" due to its
relation with the periodicity mod $N$ property.

Since (\ref{hyp12}) is a nonlinear relation, we may ask ourselves
what the ambiguity is in solving $x_3$ and $y_3$ from them. Let us
use the abbreviations
\begin{equation}
S\equiv{\sin\pi x_1\,\sin\pi x_2\over
\sin\pi y_1\,\sin\pi y_2},\qquad
T\equiv x_1+x_2-y_1-y_2,
\label{hyp13}
\end{equation}
which are single-valued functions of $x_1$, $x_2$, $y_1$, $y_2$.
We then must solve
\begin{equation}
y_3=x_3+T+2,\qquad
\sin\pi y_3=S\,\sin\pi x_3.
\label{hyp14}
\end{equation}
This has the solution
\begin{equation}
x_3={1\over2\pi\ii}\log{S-\e^{-\ii\pi T}\over S-\e^{\ii\pi T}},
\qquad y_3=x_3+T+2,
\label{hyp15}
\end{equation}
so that the only ambiguity is a translation
$x_3\mapsto x_3+M$, $y_3\mapsto y_3+M$,
shifting $x_3$ and $y_3$ by a common integer $M$. We will
use this freedom below and we conclude that we have indeed
found the required two conditions for (\ref{hyp8}) to hold.

%%%%%%%%%%%%%%%%%%%%%%%%%%%%%%%%%%%%%%%%%%%%%%%%%%%%%%%%%%%%%%%%%%%%%%%%%%%
\subsection{A double-sided hypergeometric identity\label{susec-hypd}}

Equation (\ref{hyp8}) can be further simplified after we
work out $f_2f_3/f_1$. We can do this using (\ref{hyp5})
and  $v_j-u_j=1+x_j-y_j$, for $j=1,2,3$, which follows
from (\ref{hyp7}). The result is
\begin{equation}
{f_2f_3\over f_1}={\pi^2\sin\pi u_1\over
\sin\pi x_2\sin\pi x_3\sin\pi(y_1-x_1)}
\prod_{j=1}^3{\Gamma(y_j)\Gamma(u_j)
\over\Gamma(x_j)\Gamma(v_j)\Gamma(y_j-x_j)}.
\label{hyp16}
\end{equation}
Here the functional equation of the Gamma function
$\Gamma(z)\Gamma(1-z)=\pi/\sin\pi z$ has been used.
Substituting (\ref{hyp16}) and (\ref{hyp7}) in (\ref{hyp8}),
we obtain
\begin{equation}
\sum^{\infty}_{n=-\infty}\,\prod_{j=1}^3
{\Gamma(x_j+m_j+n)\over\Gamma(y_j+m_j+n)}
={G(x_1,x_2,x_3\vert y_1,y_2,y_3)\over
\prod_{i=1}^{3}\prod_{j=1}^{3}\Gamma(y_i-x_j+m_i-m_j)},
\label{hyp17}
\end{equation}
where the function $G(x_1,x_2,x_3\vert y_1,y_2,y_3)$ is given by
\begin{equation}
G(x_1,x_2,x_3\vert y_1,y_2,y_3)\equiv{\pi^5\over
\sin\pi x_2\,\sin\pi x_3\,\prod_{i=1}^{3}\sin\pi(y_i-x_1)},
\label{hyp18}
\end{equation}
and can also be expressed as a product of ten Gamma functions.

We can simplify (\ref{hyp17}) further by absorbing the integers
$m_j$ in the $x_j$ and the $y_j$, i.e.\
$x_j+m_j\mapsto x_j$, $y_j+m_j\mapsto y_j$.
The function $G(x_1,x_2,x_3\vert y_1,y_2,y_3)$ is invariant under
this translation. We also note that in this way we utilize the
freedom in solving the two conditions (\ref{hyp10}) and
(\ref{hyp12}) in subsection \ref{susec-cond}.

As the main conclusion of this section, we find from the
large-$N$ limit of the chiral Potts model that the
following double-sided hypergeometric identity holds:
\begin{equation}
\sum_{n=-\infty}^{\infty}\prod_{i=1}^{3}
{\Gamma(x_i+n)\over\Gamma(y_i+n)}=
{G(x_1,x_2,x_3\vert y_1,y_2,y_3)\over
\prod_{i=1}^{3}\prod_{j=1}^{3}\Gamma(y_i-x_j)},
\label{dhyp}
\end{equation}
where
\begin{eqnarray}
&&G(x_1,x_2,x_3\vert y_1,y_2,y_3)=\nonumber\\
&&\qquad\prod_{j=2}^3\Gamma(x_j)\Gamma(1-x_j)
\prod_{i=1}^3\Gamma(y_i-x_1)\Gamma(1-y_i+x_1),
\label{dhyp1}
\end{eqnarray}
or equivalently (\ref{hyp18}). The identity (\ref{dhyp})
holds provided both the Saalsch\"utz condition and the
periodicity condition of subsection \ref{susec-cond} hold, i.e.
\begin{eqnarray}
&&x_1+x_2+x_3+2=y_1+y_2+y_3,\nonumber\\
&&\sin\pi x_1\,\sin\pi x_2\,\sin\pi x_3=
\sin\pi y_1\,\sin\pi y_2\,\sin\pi y_3.
\label{dhyp2}
\end{eqnarray}
Equation (\ref{dhyp}) is clearly also an identity for
\begin{equation}
\hypd{3}{x_1,x_2,x_3}{y_1,y_2,y_3}{1}=
\prod_{i=1}^{3}{\Gamma(y_i)\over\Gamma(x_i)}\,
\sum_{n=-\infty}^{\infty}\prod_{i=1}^{3}
{\Gamma(x_i+n)\over\Gamma(y_i+n)},
\label{dhyp3}
\end{equation}
which involves the double-sided hypergeometric function $\H3$
as defined e.g.\ by Slater \cite{Slater}.

We should emphasize that (\ref{dhyp}) holds whenever
(\ref{dhyp2}) is satisfied and none of the $x_j$ is an
integer. The summand and the right-hand side of (\ref{dhyp})
are meromorphic functions of their six variables and the
sum converges absolutely like $|n|^{-2}$ due to (\ref{dhyp2}),
see also subsection \ref{susec-ste1}.
From (\ref{hyp15}) we see that all solutions of (\ref{dhyp2})
are connected by analytic continuation, and the ambiguities
in solving (\ref{dhyp2}) relate to the way we go around the
logarithmic branchpoints in (\ref{hyp15}).

More generally, the $N\to\infty$ limit treated in this paper
corresponds to a $q\equiv\omega\to1$ limit of the cyclic (basic)
hypergeometric functions of \cite{AP-mf,SMS}. Looking back,
the theory of these cyclic hypergeometric functions seems to be
nearly synonymous with the theory of the integrable chiral Potts
model. The fact that we have found yet another new hypergeometric
identity confirms this point of view.

No direct proof of (\ref{dhyp}) has been given here. Such a
proof should exist and the various symmetries of (\ref{dhyp}),
\begin{itemize}
\item{Permutations of $x_1,x_2,x_3$,}
\item{Permutations of $y_1,y_2,y_3$,}
\item{Reflections $x_j\mapsto1-y_j$,
$y_j\mapsto1-x_j$, for $j=1,2,3$ simultaneously,}
\item{Translations $x_j\mapsto x_j+M$,
$y_j\mapsto y_j+M$, shifting $x_j$ and $y_j$ by a
common integer $M$, for $j=1,2$ or $3$,}
\item{Shifts $y_i\mapsto y_i+M$,
$y_j\mapsto y_j-M$, $i\ne j$, $M$ an integer,}
\end{itemize}
should be helpful in such a proof. The left-hand side and
the double product in the denominator of the right-hand side
of (\ref{dhyp}) are both separately invariant under the first
three of these symmetries, where for the reflection symmetry
it has been assumed that the conditions (\ref{dhyp2}) hold
and the functional equation of the Gamma function has to be used.
The fourth symmetry extends (\ref{dhyp}) to a star-triangle
equation, as seen above in (\ref{hyp17}), transforming one
solution of (\ref{dhyp2}) into another one. Similarly, the
fifth symmetry also maintains the validity of (\ref{dhyp}).

Therefore, the function $G(x_1,x_2,x_3\vert y_1,y_2,y_3)$
given by (\ref{hyp18}) or (\ref{dhyp1}) should also exhibit
the same symmetries, i.e.
\begin{eqnarray}
&&G(x_1,x_2,x_3\vert y_1,y_2,y_3)=
G(1\!-\!y_1,1\!-\!y_2,1\!-\!y_3\vert1\!-\!x_1,1\!-\!x_2,1\!-\!x_3)
\nonumber\\
&&\qquad=G(x_2,x_1,x_3\vert y_1,y_2,y_3)=
G(x_3,x_2,x_1\vert y_1,y_2,y_3)\nonumber\\
&&\qquad=G(x_1,x_2,x_3\vert y_2,y_1,y_3)=
G(x_1,x_2,x_3\vert y_3,y_2,y_1).
\label{dhyp4}
\end{eqnarray}
This is verified in Appendix \ref{app-D}, to which we
refer for some details.

%%%%%%%%%%%%%%%%%%%%%%%%%%%%%%%%%%%%%%%%%%%%%%%%%%%%%%%%%%%%%%%%%%%%%%%%%%%
\subsection{Relation with the Dougall--Ramanujan
formula\label{susec-doug}}

It is probably now a good idea to ask ourselves how equation
(\ref{dhyp}) fits in the existing mathematical literature. There
are not many identities available for double-sided hypergeometric
functions\footnote{We are grateful to Professor G.E.\ Andrews
for his comments on this section.} and one obvious source for them
is the textbook of Slater \cite{Slater}. It is easily seen that
the double-sided hypergeometric sum $\H{p}(\{a\};\{b\};z)$
converges only on parts of the unit circle in the $z$-plane,
and there are a few identities for $z=1$ available.

Slater gives a general $\H3$-identity following from a more
general $\H5$-identity of Dougall. The identity reads
\begin{eqnarray}
&&\hypd{3}{x_1,\phantom{1+-}x_2,\phantom{1+-}x_3}
{1+a-x_1,1+a-x_2,1+a-x_3}{1}\nonumber\\
&&\nonumber\\
&&\qquad\qquad={\Gamma(1-\half a)\Gamma(1+\half a)
\Gamma(1+\frac{3}{2}a-x_1-x_2-x_3)\over
\Gamma(1-a)\Gamma(1+a)}\nonumber\\
&&\qquad\qquad\quad\times\,
\prod_{j=1}^{3}{\Gamma(1-x_j)\Gamma(1+a-x_j)\over
\Gamma(1+a-x_1-x_2-x_3+x_j)\Gamma(1+\half a-x_j)}.
\label{dhyp5}
\end{eqnarray}
see (6.1.2.6) of \cite{Slater}, with $b=x_1$, $c=x_2$, $d=x_3$.
Note that this has four free parameters, just like in
(\ref{dhyp}) where there are two relations (\ref{dhyp2})
among the six parameters. It is easily checked that
(\ref{dhyp}) and (\ref{dhyp5}) are different, and one may say
that (\ref{dhyp5}) also has six parameters $x_j$ and $y_j=1+a-x_j$,
but with two linear relations, whereas one of the relations in
(\ref{dhyp2}) is nonlinear.

When $a=0$, (\ref{dhyp5}) reduces to the Dougall--Ramanujan
formula,
\begin{equation}
\sum_{n=-\infty}^{\infty}\prod_{i=1}^{3}
{\Gamma(x_i+n)\over\Gamma(1-x_i+n)}=
{\Gamma(1-x_1-x_2-x_3)\,\prod_{i=1}^{3}\Gamma(x_i)\over
\myunderset{1\le i<j\le3}{\prod\prod}\Gamma(1-x_i-x_j)},
\label{dhyp6}
\end{equation}
with $x_1=-x$, $x_2=-y$, $x_3=-z$ in \cite{BEd} and valid
for $\Re(x_1+x_2+x_3)<1$. Note, while comparing with
(\ref{dhyp}), that we must set $y_i=1-x_i$ and the nonlinear
periodicity condition (\ref{hyp12}) is automatically satisfied.
Imposing the Saalsch\"utz condition (\ref{hyp10}) implies
\begin{equation}
\sum_{i=1}^{3}x_i=\half,\qquad
\sum_{i=1}^{3}y_i=\textstyle{\frac{5}{2}},\qquad y_i=1-x_i,
\label{dhyp7}
\end{equation}
so that (\ref{dhyp6}) reduces to
\begin{equation}
\sum_{n=-\infty}^{\infty}\prod_{i=1}^{3}
{\Gamma(x_i+n)\over\Gamma(1-x_i+n)}=
\Gamma(\half)\,\prod_{i=1}^{3}
{\Gamma(x_i)\over\Gamma(x_i+\half)}.
\label{dhyp8}
\end{equation}
It is a simple exercise to show that the special case
$x_1+x_2+x_3=\half$, $y_i=1-x_i$ of (\ref{dhyp}) coincides with
(\ref{dhyp8}).\footnote{Only the functional relation
$\Gamma(z)\Gamma(1-z)=\pi/\sin\pi z$ and the duplication formula
$\Gamma(2z)=2^{2z-1}\Gamma(z)\Gamma(z+\half)/\Gamma(\half)$
of the Gamma function \cite{BEMOT,GRJ} are required.}

Some other interesting hypergeometric sums have been evaluated
in \cite{Searle}, now under certain nonlinear conditions, but
(\ref{dhyp}) appears still to be new.

%%%%%%%%%%%%%%%%%%%%%%%%%%%%%%%%%%%%%%%%%%%%%%%%%%%%%%%%%%%%%%%%%%%%%%%%%%%
%  Section 6
%%%%%%%%%%%%%%%%%%%%%%%%%%%%%%%%%%%%%%%%%%%%%%%%%%%%%%%%%%%%%%%%%%%%%%%%%%%
\mypr{}{\setcounter{equation}{0}}
\section{Discussion\label{sec-dis}}

There are several reasons why the above results may be of interest,
although much of this will have to be deferred to future publications.

First, the integrable chiral Potts model for finite $N$ is intimately
related with integrable deformations of a series of parafermionic
conformal field theories \cite{ZF}. It was originally proposed in
\cite{ZF} that the Fateev--Zamolodchikov model \cite{FZ} constitutes
a series of critical lattice models in the same universality classes
for each $N$. This has been numerically checked by Alcaraz
\cite{Alcaraz} for $N\le8$. The more general chiral Potts model
provides therefore lattice deformations of the conformal theory in
the chiral-field direction, see e.g.\ \cite{AM,Cardy}.

This picture should extend also to the large-$N$ limit. In section
\ref{sec-inf} we have constructed three different $N\to\infty$ limits of
the integrable chiral Potts model. These lead to three solutions of the
star-triangle equations, with an infinite sum or finite or infinite
integral, as is shown in section \ref{sec-lim}. Such exact solutions
with an infinite state space per spin ought to be of interest as they
should relate to new chiral integrable deformations of parafermionic
conformal field theories.

There are not many nontrivial exact results for nearest-neighbor
systems with infinite spin dimensionality. However, for the chiral
Potts model with finite $N$, exact results exist for the free energy
\cite{Ba,Bax-f1,Bax-f2,Bax-f3}, order parameter \cite{AlMPT},
groundstate energy and excitation spectra \cite{AlMP,Bax-su,MR}
of the associated quantum chain, and surface tensions
\cite{Bax-if1,Bax-if2,AP-rev,OB}. Many of these results have been
obtained using systems of functional equations for transfer matrices
\cite{BS,BBP}.

Large-$N$ limits of these quantities can be constructed. For
example, we can use a conjecture for the order parameters \cite{AlMPT},
i.e.\ for each $n$ with $1\le n\le N-1$ an order parameter is given
by the expectation value 
\begin{equation}
\langle\omega^{na_0}\rangle={(1-{k'}^2)}^{n(N-n)/2N^2},
\label{conc1}
\end{equation}
where $a_0=0,\cdots,N-1$ is the random value of a given bulk spin,
say at the origin. In the limiting regime I, where this spin $a_0$
can run through all positive and negative integers, (\ref{conc1})
tends to
\begin{equation}
\langle\e^{\ii x a_0}\rangle={(1-{k'}^2)}^{x(2\pi-x)/8\pi^2},
\label{conc2}
\end{equation}
where $x$ is a real number $0<x<2\pi$ resulting from the limit
$2\pi n/N\to x$. Similar limits can be constructed from existing
exact results for some other thermodynamic quantities.

The precise status of (\ref{conc2}) must still be determined as
the above construction presumes the interchange of the $N\to\infty$
limit with the thermodynamic limit followed by the field limit
defining the order parameter. Furthermore, it would be interesting
to study all these thermodynamic quantities further within a larger
$\infty$-state model containing the integrable manifold and we hope
to return to this in the future.

Finally, from a mathematical point of view, the $N\to\infty$ limit
corresponds to a $q\equiv\omega\to1$ limit of the cyclic (basic)
hypergeometric functions of \cite{AP-mf,SMS}, which are intimately
related with the integrable chiral Potts model. However, the
connection of the double-sided series (\ref{infi27}) and (\ref{dhyp})
with more general cyclic hypergeometric series will be treated elsewhere.

%%%%%%%%%%%%%%%%%%%%%%%%%%%%%%%%%%%%%%%%%%%%%%%%%%%%%%%%%%%%%%%%%%%%%%%%%%%
%  Acknowledgments
%%%%%%%%%%%%%%%%%%%%%%%%%%%%%%%%%%%%%%%%%%%%%%%%%%%%%%%%%%%%%%%%%%%%%%%%%%%
\section*{Acknowledgments}

We thank our hosts Professors A.J.\ Guttmann, P.A.\ Pearce,
and P.J.\ Forrester at the University of Melbourne and
Professors H.W.\ Capel and B.\ Nienhuis at the University
of Amsterdam for much support. We thank Professor G.E.\ Andrews
for helpful discussions and our many other colleagues in Melbourne
and Amsterdam for their warm hospitality.
This work has been supported in part by NSF Grants
Nos.\ PHY 97--22159 and PHY 97--24788. Some further financial
support from The University of Melbourne and from FOM
(`Stichting voor Fundamenteel Onderzoek der Materie'),
which is financially supported by the NWO (`Nederlandse Organisatie
voor Wetenschappelijk Onderzoek'), is gratefully acknowledged.

%%%%%%%%%%%%%%%%%%%%%%%%%%%%%%%%%%%%%%%%%%%%%%%%%%%%%%%%%%%%%%%%%%%%%%%%%%%
%  Section A
%%%%%%%%%%%%%%%%%%%%%%%%%%%%%%%%%%%%%%%%%%%%%%%%%%%%%%%%%%%%%%%%%%%%%%%%%%%
\mypr{}{\setcounter{equation}{0}}
\appendix
\section{The Constant {\em A}\label{app-A}}

In this appendix we verify that the formula for constant $A$
given in (\ref{infi1}) works out in all four cases.

Substituting (\ref{infi2}) and (\ref{infi3}) we get, using
(\ref{repa3}) and (\ref{repa8}),
\begin{eqnarray}
(A_{\pq})^2&=&
{\cos^2\half(\theta_{\p}-\phi_{\q})\over
\cos^2\half(\theta_{\q}-\phi_{\p})}=
{1+\cos\theta_{\p}\cos\phi_{\q}+\sin\theta_{\p}\sin\phi_{\q}
\over
1+\cos\theta_{\q}\cos\phi_{\p}+\sin\theta_{\q}\sin\phi_{\p}}
\nonumber\\&=&
{1+k^2+2k\cos\theta_{\p}\over1+k^2+2k\cos\theta_{\q}}
={\sin\theta_{\p}\sin\phi_{\q}\over
\sin\phi_{\p}\sin\theta_{\q}}=
\e^{2\gamma_{\p}-2\gamma_{\q}},
\label{infi6}
\end{eqnarray}
and
\begin{eqnarray}
(\bA_{\pq})^2&=&
{\sin^2\half(\theta_{\q}-\theta_{\p})\over
\sin^2\half(\phi_{\q}-\phi_{\p})}=
{1-\cos\theta_{\p}\cos\theta_{\q}-\sin\theta_{\p}\sin\theta_{\q}
\over
1-\cos\phi_{\q}\cos\phi_{\p}-\sin\phi_{\q}\sin\phi_{\p}}
\nonumber\\&=&
{(1+k^2+2k\cos\theta_{\p})(1+k^2+2k\cos\theta_{\q})\over1-k^2}
\nonumber\\&=&{\sin\theta_{\p}\sin\theta_{\q}\over
\sin\phi_{\p}\sin\phi_{\q}}=
\e^{2\gamma_{\p}+2\gamma_{\q}},
\label{infi7}
\end{eqnarray}
so that we receive agreement with (\ref{repa11}) and (\ref{repa12}).
The signs of the square roots are easily verified for the limit
$\theta_{\q}\to\theta_{\p}$ and $\phi_{\q}\to\phi_{\p}$
in case of (\ref{infi6}), or
$\phi_{\p}\to\theta_{\p}$ and $\phi_{\q}\to\theta_{\q}$ for
case (\ref{infi7}).

Similarly, substituting (\ref{infi8}) and (\ref{infi9}) into
the constants $A$ of (\ref{infi1}) and using (\ref{repa7}) we have
\begin{eqnarray}
\Af_{\pq}&=&{\sin\half(\tb_{\p}-\tb_{\q})\over
\sin\half(\pb_{\p}-\pb_{\q})}=
\e^{\ii\pi(\lambda_{\p}+\lambda_{\q})}
{\e^{\gamma_{\q}-\ii\pi\lambda_{\q}}-
\e^{\gamma_{\p}-\ii\pi\lambda_{\p}}\over
\e^{\gamma_{\q}+\ii\pi\lambda_{\q}}-
\e^{\gamma_{\p}+\ii\pi\lambda_{\p}}}\nonumber\\&=&
\e^{\ii(\phi_{\p}+\phi_{\q}-\theta_{\p}-\theta_{\q})/2},
\label{infi10}
\end{eqnarray}
\begin{eqnarray}
\bAf_{\pq}&=&{\cos\half(\pb_{\p}-\tb_{\q})\over
\cos\half(\tb_{\p}-\pb_{\q})}=
\e^{\ii\pi(\lambda_{\q}-\lambda_{\p})}
{\e^{\gamma_{\p}+\ii\pi\lambda_{\p}}
\e^{\gamma_{\q}-\ii\pi\lambda_{\q}}+1\over
\e^{\gamma_{\p}-\ii\pi\lambda_{\p}}
\e^{\gamma_{\q}+\ii\pi\lambda_{\q}}+1}\nonumber\\&=&
\e^{\ii(\phi_{\q}-\theta_{\q}-\phi_{\p}+\theta_{\p})/2},
\label{infi11}
\end{eqnarray}
in agreement with (\ref{repa13}) and (\ref{repa14}).

%%%%%%%%%%%%%%%%%%%%%%%%%%%%%%%%%%%%%%%%%%%%%%%%%%%%%%%%%%%%%%%%%%%%%%%%%%%
%  Section B
%%%%%%%%%%%%%%%%%%%%%%%%%%%%%%%%%%%%%%%%%%%%%%%%%%%%%%%%%%%%%%%%%%%%%%%%%%%
\mypr{}{\setcounter{equation}{0}}
\section{Mathematical Details of Large-{\em N} Limit\label{app-N}}

In this appendix we prove that
\begin{equation}
S_n(\alpha)\equiv\sum_{j=1}^{n}\log{\sin[\pi(j\!+\!\alpha\!-\!1)/N]\over
\pi(j\!+\!\alpha\!-\!1)/N}
\label{A-sum}
\end{equation}
is asymptotically for large $N$ given by
\begin{equation}
S_n(\alpha)\approx\sum_{j=0}^{\infty}{\B_j(\alpha)\over j!}
{\left({\pi\over N}\right)}^{j-1}
\left[\phi^{(j)}(\pi n/N)-\phi^{(j)}(0)\right],
\label{A-asym}
\end{equation}
where
\begin{equation}
\phi(z)\equiv\int_{0}^{z}\d z'\log{\sin z'\over z'}
\label{B1}
\end{equation}
and $\phi^{(j)}(z)$ is its $j$th derivative.

The expression (\ref{A-sum}), valid for $0\le n<N$, must be
replaced by
\begin{equation}
S_n(\alpha)=-\sum_{j=1}^{-n}\log{\sin[\pi(j\!-\!\alpha)/N]\over
\pi(j\!-\!\alpha)/N}=-S_{N-n}(1\!-\!\alpha),
\label{B2}
\end{equation}
for $-N<n\le0$, which has the identical asymptotic expansion
(\ref{A-asym}). This can be easily verified using
\begin{equation}
\B_j(1-\alpha)=(-1)^{j}\B_j(\alpha)
\label{B3}
\end{equation}
and
\begin{equation}
\phi^{(j)}(-z)-\phi^{(j)}(0)=(-1)^{j+1}[\phi^{(j)}(z)-\phi^{(j)}(0)].
\label{B4}
\end{equation}

From elementary calculus we note that the derivatives
\begin{equation}
{\d\phi(z)\over\d z}=\log{\sin z\over z},\qquad
{\d^2\phi(z)\over\d z^2}=\cot z-{1\over z},
\label{B5}
\end{equation}
have Taylor expansions involving the Bernoulli numbers
\cite{BEMOT,GRJ,GRJa}. We have
\begin{eqnarray}
&\phi(z)&=
\sum_{l=2}^{\infty}\,{(2i)^l\B_{l}\,z^{l+1}\over l\,(l\!+\!1)!}=
\sum_{k=1}^{\infty}\,{(-1)^k2^{2k}\B_{2k}\,z^{2k+1}\over
2k\,(2k\!+\!1)!}\nonumber\\&&=
-\sum_{k=1}^{\infty}\,{\zeta(2k)z^{2k+1}\over k(2k\!+\!1)\pi^{2k}},
\label{A-phi}
\end{eqnarray}
convergent for $\vert z\vert<\pi$. Here $\zeta(s)$ is the Riemann
zeta function.

Therefore, we can write
\begin{equation}
S_n(\alpha)=\sum_{j=1}^{n}\sum_{l=2}^{\infty}\,
{\left({2\pi\ii\over N}\right)}^l\,{\B_{l}\over l\!\cdot\!l!}\,
(j\!+\!\alpha\!-\!1)^l.
\label{B6}
\end{equation}
For
\begin{equation}
-N<\alpha,n+\alpha-1<N,
\label{A-conv}
\end{equation}
this is absolutely convergent, allowing
us to perform the sum over $l$ in terms of Bernoulli polynomials
\cite{BEMOT,GRJ},
\begin{equation}
\B_{n}(x)\equiv{}^{\hbox{``}}{(x+\B)}^n\,{}^{\hbox{"}}\equiv
\sum_{k=0}^{n}\,{n\atopwithdelims()k}\,
\B_{k}\,x^{n-k},
\label{B7}
\end{equation}
with the properties
\begin{equation}
\B_{n}(x+1)=\B_{n}(x)+nx^{n-1},
\label{B8}
\end{equation}
\begin{equation}
\sum_{j=1}^{n}\,(j\!+\!\alpha\!-\!1)^l=
{\B_{l+1}(\alpha\!+\!n)-\B_{l+1}(\alpha)\over l\!+\!1},
\label{B9}
\end{equation}
\begin{equation}
\B_{l+1}(\alpha\!+\!n)=
{}^{\hbox{``}}{(\alpha\!+\!n\!+\!\B)}^{l+1}\,{}^{\hbox{"}}=
\sum_{k=0}^{l+1}\,{l\!+\!1\atopwithdelims()k}\,
\B_{k}(\alpha)\,x^{l+1-k}.
\label{B10}
\end{equation}
We find therefore
\begin{equation}
S_n(\alpha)=\sum_{l=2}^{\infty}\sum_{k=0}^{l}\,
{\left({2\pi\ii\over N}\right)}^l\,{\B_{l}\over l}\,
{\B_{k}(\alpha)\over k!}\,\left({\d\over\d x}\right)^k
\left.{x^{l+1}\over(l\!+\!1)!}\right|_{x=n}.
\label{B11}
\end{equation}
Here the sum over $l$ converges absolutely within range (\ref{A-conv}).

However, the double sum converges only relatively, as one can also
verify numerically. In order to obtain the result (\ref{A-asym}) we
need to interchange the two sums, leading to an asymptotic expansion
as a consequence. We find
\begin{equation}
S_n(\alpha)=\sum_{k=0}^{\infty}\,
{\B_{k}(\alpha)\over k!}\,{N\over\pi}\,
\left({\d\over\d x}\right)^k\,
\sum_{\textstyle{l=k\atop l\ge2}}^{\infty}\,
{(2\ii)^l\,\B_{l}\over l\!\cdot\!(l\!+\!1)!}
\left.{\left({\pi x\over N}\right)}^{l+1}\right|_{x=n},
\label{B12}
\end{equation}
from which (\ref{A-asym}) immediately follows, or equivalently
\begin{eqnarray}
S_n(\alpha)&\approx&{N\over\pi}\,\phi(\pi n/N)
+(\alpha-\half)\phi^{(1)}(\pi n/N)\nonumber\\
&&+\sum_{l=1}^{\infty}{\B_{l+1}(\alpha)\over(l+1)!}
{\left({\pi\over N}\right)}^{l}
\left[\phi^{(l+1)}(\pi n/N)-\phi^{(l+1)}(0)\right],
\label{A-asym2}
\end{eqnarray}
after using $\B_0(\alpha)=1$, $\B_1(\alpha)=\alpha-{1\over2}$,
$\phi(0)=\phi^{(1)}(0)=0$.

The first line of (\ref{A-asym2}) is sufficient for our purposes, as
we need to keep only terms of order $\O(1)$ as $N\to\infty$.
Therefore, we conclude this appendix with estimating the second line.
For $n\ge2$, $0\le x\le1$, the Bernoulli polynomial can be
expressed as \cite{BEMOT,GRJ}
\begin{equation}
\B_n(x)=-{2\!\cdot\!n!\over(2\pi)^n}\,
\sum_{k=1}^{\infty}\,{\cos(2\pi k x\!-\!\half\pi n)\over k^n},
\label{B13}
\end{equation}
so that
\begin{equation}
\vert\B_n(x)\vert\le{2\!\cdot\!n!\over(2\pi)^n}\,\zeta(n),
\quad\hbox{for }n\ge2,\quad0\le x\le1,
\label{A-bound1}
\end{equation}
with $1<\zeta(n)\le\zeta(2)=\pi^2/6$. Also, for $m=2,3,\cdots,$ and
$\vert z\vert<\pi$, we have from (\ref{A-phi})
\begin{equation}
\phi^{(m)}(z)-\phi^{(m)}(0)=
-\sum_{k=
%\left\lceil\vphantom{\sum}\right.\half m
%\left.\vphantom{\sum}\right\rceil
\lceil m/2\rceil}^{\infty}\,
{2\!\cdot\!(2k\!-\!1)!\,\zeta(2k)z^{2k+1-m}
\over(2k\!+\!1\!-\!m)!\,\pi^{2k}},
\label{A-phi2}
\end{equation}
where $\lceil x\rceil$ is the smallest integer $\ge x$.
Since all the coefficients in (\ref{A-phi2}) are negative, we have
\begin{equation}
\phi^{(m)}(z)-\phi^{(m)}(0)=(-1)^{m+1}[\phi^{(m)}(-z)-\phi^{(m)}(0)]\le0
\label{B14}
\end{equation}
and monotonically decreasing for $0\le z\le\pi$. Therefore,
we can estimate (\ref{A-phi2}) by replacing all $\zeta(2k)$ by
its minimum 1 or its maximum $\zeta(2\lceil\half m\rceil)$,
which makes the sum a binomial-type expansion. We find
\begin{equation}
\chi^{(m)}(z)\le
-({\sign}\,z)^{m+1}\left[\phi^{(m)}(z)-\phi^{(m)}(0)\right]\le
\zeta(2\lceil\half m\rceil)\,\chi^{(m)}(z),
\label{A-bound2}
\end{equation}
where
\begin{eqnarray}
\chi^{(2n)}(z)&\equiv&(2n-2)!\,\left[{1\over(\pi-\vert z\vert)^{2n-1}}-
{1\over(\pi+\vert z\vert)^{2n-1}}\right],\nonumber\\
\chi^{(2n+1)}(z)&\equiv&(2n-1)!\,\left[{1\over(\pi-\vert z\vert)^{2n}}+
{1\over(\pi+\vert z\vert)^{2n}}-{2\over\pi^{2n}}\right],
\label{B15}
\end{eqnarray}
which are both positive for $n=1,2,\cdots$.

Bounds (\ref{A-bound1}) and (\ref{A-bound2}) are quite sharp, as can
also be numerically verified. The bounds show that (\ref{A-asym}) is
an asymptotic expansion, with absolute value of terms roughly bounded
by $l!/(\pi N)^l$ and error less than the bound on the first ignored
term. This is true for $N\ge2$ and $\vert n\vert<N$. But the
expansion becomes particularly useful for large $N$ and the full range
of $n$ taken as $\vert n\vert\le\half N$.

%%%%%%%%%%%%%%%%%%%%%%%%%%%%%%%%%%%%%%%%%%%%%%%%%%%%%%%%%%%%%%%%%%%%%%%%%%%
%  Section C
%%%%%%%%%%%%%%%%%%%%%%%%%%%%%%%%%%%%%%%%%%%%%%%%%%%%%%%%%%%%%%%%%%%%%%%%%%%
\mypr{}{\setcounter{equation}{0}}
\section{Derivation of (\ref{lim8})\label{app-F}}

The factor $F\equiv F_{\pq}$ is obtained substituting
(\ref{lim4})--(\ref{lim8}) into (\ref{lim3}). We can simplify it
using (\ref{lim1}) and the functional equation of the Gamma function
$\Gamma(x)\Gamma(1-x)=\pi/\sin\pi x$. This gives
\begin{equation}
F=2\,{\biggl({N\over2\pi}\biggr)}^{\beta-\alpha}\,
{\bA^{1/2}\,\mbox{$\bAf$}^{1/2}\over A^{1/2}}\,
{\Gamma(\alpha)\Gamma(\bbeta)\Gamma(1-\bbeta+\balpha)\sin\pi\balphaf
\over\Gamma(\balpha)\Gamma(\beta)}.
\label{C1}
\end{equation}
Part of this expression can be simplified further using
$A\equiv\sin\pi\beta/\sin\pi\alpha$, (\ref{lim8}) and (\ref{lim1}).
More precisely,
\begin{eqnarray}
4\bA\,\bAf\,\sin^2\pi\balphaf&&=
4\bA\,\sin\pi\balphaf\sin\pi\bbetaf\nonumber\\
&&\nonumber\\
&&=
2\bA\cos\pi(\balphaf-\bbetaf)-2\bA\cos\pi(\balphaf+\bbetaf)\nonumber\\
&&\nonumber\\
&&=-2\bA\cos\pi(\balpha-\bbeta)+2\bA\cosh\log\bA\nonumber\\
&&\nonumber\\
&&=\bA^2-2\bA\cos\pi(\balpha-\bbeta)+1=
{\sin^2\pi(\bbeta-\balpha)\over\sin^2\pi\balpha}.
\label{C2}
\end{eqnarray}
Taking the square root of this in domain (\ref{lim2}) and applying
it to (\ref{C1}) we find
\begin{equation}
F={\biggl({N\over2\pi}\biggr)}^{\beta-\alpha}\,
{\Gamma(\alpha)\Gamma(\bbeta)\Gamma(1-\bbeta+\balpha)
\sin\pi(\bbeta-\balpha)
\over A^{1/2}\Gamma(\balpha)\Gamma(\beta)\sin\pi\balpha}.
\label{C3}
\end{equation}
This can be further reduced to (\ref{lim8}) using the functional
equation of the Gamma function again.

%%%%%%%%%%%%%%%%%%%%%%%%%%%%%%%%%%%%%%%%%%%%%%%%%%%%%%%%%%%%%%%%%%%%%%%%%%%
%  Section D
%%%%%%%%%%%%%%%%%%%%%%%%%%%%%%%%%%%%%%%%%%%%%%%%%%%%%%%%%%%%%%%%%%%%%%%%%%%
\mypr{}{\setcounter{equation}{0}}
\section{Verification of (\ref{dhyp4})\label{app-D}}

In this appendix we show that under the conditions (\ref{dhyp2}), or
\begin{equation}
2+\sum_{j=1}^{3}x_j=\sum_{j=1}^{3}y_j,\quad
\prod_{j=1}^{3}\sin\pi x_j=\prod_{j=1}^{3}\sin\pi y_j
\equiv\tau/(2\ii)^3.
\label{D0}
\end{equation}
the symmetry relations (\ref{dhyp4}) hold or, equivalently,
\begin{equation}
{\sin\pi(y_1-x_1)\sin\pi(y_2-x_1)\sin\pi(y_3-x_1)
\over\sin\pi x_1}\equiv\sigma/(2\ii)^2
\label{D1}
\end{equation}
is fully symmetric both in $\{x_1,x_2,x_3\}$ and in $\{y_1,y_2,y_3\}$
and it is also invariant under $x_j\mapsto1-y_j$,
$y_j\mapsto1-x_j$, for $j=1,2,3$ simultaneously.
Indeed this gives (\ref{dhyp4}), as $G(x_1,x_2,x_3\vert y_1,y_2,y_3)=
(2\pi\ii)^5\tau/\sigma$.

With the definitions
\begin{equation}
\xi_i\equiv\e^{\ii\pi x_i},\quad\eta_i\equiv\e^{\ii\pi y_i},
\quad\hbox{for }i=1,2,3,
\label{D2}
\end{equation}
the conditions (\ref{D0}) become
\begin{equation}
\xi_1\xi_2\xi_3=\eta_1\eta_2\eta_3\equiv\rho,\quad
\prod_{i=1}^{3}(\xi_i-\xi_i^{-1})=
\prod_{i=1}^{3}(\eta_i-\eta_i^{-1})=\tau.
\label{D3}
\end{equation}
We can now expand $\rho\tau$ and rearrange terms. We obtain
\begin{equation}
\rho^2\sum_{i=1}^{3}(\xi_i^{-2}-\eta_i^{-2})=
\sum_{i=1}^{3}(\xi_i^2-\eta_i^2).
\label{D4}
\end{equation}
We use this to expand $\sigma$. Successively, we find
\begin{eqnarray}
\sigma&&={\rho^2-\xi_1^2\rho^2\sum_{i=1}^{3}\eta_i^{-2}
+\xi_1^4\sum_{i=1}^{3}\eta_i^2-\xi_1^6\over
\rho\xi_1^2(\xi_1^2-1)}\nonumber\\
&&\nonumber\\
&&={\xi_2^2\xi_3^2-\rho^2\sum_{i=1}^{3}\xi_i^{-2}
+\sum_{i=1}^{3}\xi_i^2-\sum_{i=1}^{3}\eta_i^2
+\xi_1^2\sum_{i=1}^{3}\eta_i^2-\xi_1^4\over
\rho(\xi_1^2-1)}\nonumber\\
&&\nonumber\\
&&=\sum_{i=1}^{3}(\eta_i^2-\xi_i^2)/\rho.
\label{D5}
\end{eqnarray}
This shows that $\sigma$ has the required permutation symmetries.
The invariance of $\sigma$ under the reflection symmetry
\begin{equation}
\xi_i\mapsto-\eta_i^{-1},\quad
\eta_i\mapsto-\xi_i^{-1},\quad
(i=1,2,3),\quad\hbox{and }
\rho\mapsto-\rho^{-1},
\label{D6}
\end{equation}
follows from (\ref{D4}) and (\ref{D5}). As $\tau$ obeys these
symmetries trivially, we can now complete the proof of (\ref{dhyp4}).

%%%%%%%%%%%%%%%%%%%%%%%%%%%%%%%%%%%%%%%%%%%%%%%%%%%%%%%%%%%%%%%%%%%%%%%%%%%
%  Section R
%%%%%%%%%%%%%%%%%%%%%%%%%%%%%%%%%%%%%%%%%%%%%%%%%%%%%%%%%%%%%%%%%%%%%%%%%%%
%\nonumsection{References}

\end{document}